\documentclass[11pt]{article}   % For Latex2e
\usepackage{amssymb,amscd,latexsym}   % For Latex2e
\usepackage{amsmath}
\usepackage{amsthm}
\newcommand{\rar}{\rightarrow}
\newcommand{\lar}{\longrightarrow}
\newcommand{\llar}{-\kern-5pt-\kern-5pt\longrightarrow}
\newcommand{\surjects}{\twoheadrightarrow}

\newtheorem{Theorem}{Theorem}[section]
\newtheorem{Lemma}[Theorem]{Lemma}
\newtheorem{Corollary}[Theorem]{Corollary}
\newtheorem{Proposition}[Theorem]{Proposition}
\newtheorem{Remark}[Theorem]{Remark}
\newtheorem{Example}[Theorem]{Example}

\newtheorem{Definition}[Theorem]{Definition}
\newtheorem{Question}[Theorem]{Question}
\def\sqr#1#2{{\vcenter{\hrule height.#2pt
        \hbox{\vrule width.#2pt height#1pt \kern#1pt
            \vrule width.#2pt}
        \hrule height.#2pt}}}
\def\phi{\varphi}
\def\demo{\noindent{\bf Proof. }}
\def\square{\mathchoice\sqr64\sqr64\sqr{4}3\sqr{3}3}
\def\qed{\hspace*{\fill} $\square$}

\def\floor#1{{\left \lfloor{#1}\right \rfloor}}

\DeclareMathOperator{\Hom}{Hom}

\DeclareMathOperator{\Ass}{Ass}

\DeclareMathOperator{\Ext}{Ext}

\DeclareMathOperator{\depth}{depth}

\def\xx{{\bf x}}
\def\yy{{\bf y}}

\def\XX{{\bf X}}
\def\YY{{\bf Y}}

\def\ZZ{{\bf Z}}
\def\ff{{\bf f}}

\def\ff{{\bf f}}
\def\gg{{\bf g}}

\def\fm{{\mathfrak m}}

\def\hht{{\rm ht}\,}
\def\depth{{\rm depth}\,}

\def\ker{{\rm ker}\,}

\def\Ext#1#2#3#4{{\rm Ext}\,^{#1}_{#2}({#3},{#4})}

\def\restr{{\kern-1pt\restriction\kern-1pt}}

\def\pp{{\mathbb P}}

 2
 3
 4
 5
 4
 3
 2
 1
 1
 2
 3
 1
 0
 1
 2
 3
 4
 5

\begin{document}
\begin{center}
{\Large{\bf\sc Symbolic powers of perfect ideals of codimension $2$}}\\[3pt]
{\Large{\bf\sc and birational maps }}
\footnotetext{2010 AMS {\it Mathematics Subject
Classification}: 13A30, 13C14, 13C40, 13D02, 13E15, 13H10, 13P10, 14E05, 14E07.}

\vspace{0.3in}

{\large\sc Zaqueu Ramos}\footnote{Parts of this work were done while this author held a
Doctoral Fellowship (CAPES, Brazil).} \quad\quad
 {\large\sc Aron  Simis}\footnote{Partially
supported by a CNPq grant.}

\end{center}

\begin{center}

{\em To Wolmer Vasconcelos on his 75th birthday, for his seminal mathematical ideas}
\end{center}

%\tableofcontents

\bigskip

\begin{abstract}

This work is about symbolic powers of codimension two perfect ideals
in a standard polynomial ring over a field, where the entries of the corresponding
presentation matrix are general linear forms.
The main contribution of the present approach
is the use of the birational theory underlying
the nature of the ideal and the details of a deep interlacing between generators of its symbolic powers and the
inversion factors stemming from the inverse map to the birational map defined by the linear system spanned by
the generators of this ideal.
A full description of the corresponding symbolic Rees algebra is given in some cases.
An application is an affirmative solution of a conjecture of Eisenbud--Mazur in \cite[Section 2]{EiMa}.

\end{abstract}

\section*{Introduction}

Let $I\subset R$ denote an ideal in a Noetherian ring and let $r\geq 0$ be an integer.
The  $r$th symbolic power $I^{(r)}$ of $I$
can be defined as the inverse image of $S^{-1}I^{r}$ under the natural homomorphism $R\rar S^{-1}R$
of fractions, where $S$ is the complementary set of the union of the associated primes of $R/I$.
There is a known hesitation as to whether one should take the whole set of associated primes of $R/I$ or just
its minimal primes or even those of minimal codimension or maximal dimension.
In this work we need not worry about this dilemma because the notion will only be employed in the case
of a codimension $2$ perfect ideal in a Cohen--Macaulay ring -- actually, a polynomial ring over a field.
In this setup there is no ambiguity and $I^{(r)}$ is precisely the intersection of the primary components
of the ordinary power $I^{r}$ relative to the associated primes of $R/I$, i.e., the unmixed part of $I^{r}$.

A more serious problem is the characteristic of the base field.
In characteristic zero, if $I$ is a radical ideal, one has the celebrated Zariski--Nagata differential
characterization of $I^{(r)}$ (see  \cite[3.9]{E} and the references there).
The subject in positive characteristic or mixed characteristic gives a quite different panorama, often much harder but with
different methods anyway.
Essential parts of this work assume characteristic zero.
This is not due to a need of using the Zariski--Nagata criterion upfront, but rather to an urge
of dealing with Jacobian matrices and
using Bertini's theorem.
Many technical results will be valid just over an infinite field, hence there has been an effort
to convey when the characteristic is an issue at specific places.
On the other hand, since we will draw quite substantially on aspects of birational maps, it may be a good idea in
those instances to think about $k$ as being algebraically closed.

The main object of concern is an $m\times (m-1)$ matrix whose entries are general $1$-forms
in a polynomial ring $R=k[X_1,\ldots, X_n]$ over an infinite field $k$ -- called herein {\em general linear matrices}.
We will focus on the ideal $I\subset R$
generated by the $(m-1)$-minors of the matrix.
The group ${\rm Gl}(m,k)\times {\rm Gl}(n,k)\times {\rm Gl}(m-1,k)$ acts on the set of all linear
$m\times (m-1)$ matrices over $k$.
Depending on the relative values of $m$ and $n$, these matrices may fail to be
 $1$-generic in the sense of \cite[Proposition-Definition 1.1]{Eisenbud2}.
 For $m\times (m-1)$ linear matrices the condition of being $1$-generic implies that $n\geq 2(m-1)$.
 Though natural in various contexts, $1$-genericity goes only ``half'' way  the cases.
 For $n < m$, e.g., the above triple action does not preserve the property of being general linear, as is clear
 that one may introduce a certain number of zero entries in the resulting matrix.
 The property is preserved if the action is restricted to a suitable open set
 of ${\rm Gl}(m,k)\times {\rm Gl}(n,k)\times {\rm Gl}(m-1,k)$.
 As a simple example, take $m=2, n=1$. Then the $2\times 1$ matrix $(\alpha x, \beta x)^t\, (\alpha\neq 0)$
 can be converted to $(\alpha x,0)^t$ by the left action of the element
 $$\mathfrak{g}=\left(\begin{matrix}
 1 & 0\\
 -\beta/\alpha & 1
 \end{matrix}
 \right)
 $$
 (identified with $\mathfrak{g}\times \boldsymbol 1\times \boldsymbol 1$),
 but not if the acting group element has general coefficients.
This scrambling in the orbit of a general linear matrix in the present sense is often a root of difficulty when handling
ideal theoretic properties stemming from the data.

For practical purposes, a set of general $1$-forms of the sort we assume can be taken to be a set of random $k$-linear
combinations of the variables.
Any such set of cardinality $m(m-1)$ can be ordered as the entries of a matrix,
so there are plenty of such matrices.
However, for the sake of subsequent ideal theoretic development we introduce a more formal
definition right at the outset (see Definition~\ref{DEF}).
A weaker form would require that the entries individually be general forms in general linear position
(i.e., every subset of the entries of cardinality at most $n$
be  $k$-linearly independent).
There are even weaker forms that have been considered in the literature.
It is not clear whether all the main results of the paper go through in those environments.
Examples are given to show that some of the crucial results obtained no longer subsist
in a less general frame.

The idea behind the present subject is akin to other places where one introduces an object in terms
of suitable general data -- such is the case of the notion of the generic initial ideal of a polynomial ideal in
Gr\"obner base theory (see, e.g., \cite{Conca_et_al}).
One starts out with some random
like definition and then pursues some well-defined algebraic behavior for these data.
If one thinks about it, the apparent difficulty surfaces at once.
This often justifies why some of the arguments spelled in such a setup are long and detailed, whereas they often appear
nearly obvious.

\medskip

Now, for a homogeneous ideal $I\subset R$ generated in fixed degree,  whose syzygies are generated by
``enough'' linear syzygies, its generators are very close to span a linear system defining a
birational map from a projective space onto its image.
This strategy has been largely explored in recent years by several authors.
Thus, the details of the geometry of birational maps can be accommodated in terms of numerical invariants
from commutative algebra.
However,  finding room in this accommodation for symbolic powers has not, to our knowledge, been brought up so far.
This is one of our main observations in this work. Together with a good grip of the algebraic and homological
properties of the base ideal $I$, it constitutes the main bulk of the paper.

The main results of this paper are shown in Theorem~\ref{specialization}, Theorem~\ref{associated_primes},
Theorem~\ref{symbolic_m=n}, Theorem~\ref{funny_module}, Theorem~\ref{The D's} and Theorem~\ref{symbolic_m=n+1}.
A consequence of Theorem~\ref{symbolic_m=n} is a solution,  over a field of characteristic zero,
of a conjecture stated by Eisenbud--Mazur in \cite[Section 2]{EiMa}  -- we are not aware of a previous solution
in the literature.

\smallskip

We now briefly describe the contents of each section.

The first section, divided in two parts, gives an overview of the basic material on symbolic powers and
on birational maps.
The first part gives the tool used to approach the nature of
the symbolic algebra in the present context. It is based on  an idea of Vasconcelos that brought in
the management of the ideal transform in this setup.
The second part discusses a couple of useful facts, apparently thus far unnoticed in such generality.
They have to do with the so-called
{\em inversion factor} of a birational map. These properties are proved in
Lemma~\ref{uniqueness_of_target_factor}, Proposition~\ref{jac_vs_factor} and 
Proposition~\ref{inversionfactor_is_symbolic}.

The second section contains the main  results of the paper.
It starts with some preliminaries on a perfect ideal $I\subset R=k[X_1,\ldots, X_n] (n\geq 3)$ of codimension $2$
whose structural  $m\times (m-1)$ matrix is a general linear matrix.
We first show that the other Fitting ideals attain an expected codimension and that $I$ enjoys  typical properties
which depend on the values of $m,n$.
Thus, for $n\geq 4$, $I$ is a normal prime ideal provided char$(k)=0$ (and possibly in general); moreover, it is of linear type
if (and only if) $m\leq n$ and it is normally torsionfree if and only if $m<n$.
Therefore, such an ideal is only really of new substance in the case where $m\geq n$.
In the sequel we show that $I$ satisfies a generalized property of Artin--Nagata, called $(G_n)$
and that, for any exponent $r$, the symbolic power $I^{(r)}$ coincides with the $(\XX)$-saturation of $I^r$
(in other words, the unmixed part of $I^r$ is its saturation).
Since $I$ is prime for $n\geq 4$, the symbolic power  $I^{(r)}$  is just the $I$-primary component of $I^r$
and the latter has at most one further  associated prime, namely, $(\XX)$.
To go one step forward, we introduce certain graded pieces of the approximation complex, along with other techniques and
a recent result of A. B. Tchernev, to deduce that if $I^{(r)}\neq I^r$ then necessarily $r\geq n-1$.
This result becomes an important tool for the rest of the work.

So much for the main ideal theoretical and homological properties.
On a second part of the same section, we deal with the `equations' of $I$.
Namely, we bring up the role of Rees algebra of $I$ in the underlying birational map
based on the linear system spanned by the generators of $I$.
Specifically, we show that for $m\geq n\geq 3$ the ideal $I$ is the base ideal of a birational map of $\pp^{n-1}$ onto
the image in $\pp^{m-1}$.
This result is based on the {\em fiber type} nature of $I$ -- i.e., its Rees algebra is simplest beyond the linear type
situation -- and on a special case of the criterion of birationality established in \cite{AHA}.

In this part we bring in detail the role of the inversion factors associated to the birational map in question, showing that
they are natural elements in the symbolic power $I^{(n-1)}$ not belonging to the ordinary power $I^{n-1}$.
Inversion factors have appeared before in the classical theory of plane Cremona maps, where they are a version
of the {\em principal curves} (see, e.g., \cite[Chapter 3]{alberich}).
However, to our knowledge the notion has never been explicitly addressed for Cremona maps in higher dimension,
much less for birational maps onto their images (classically called ``rational representations'' of projective space).
We introduce them here in this larger generality and dimension. A bit surprisingly, they keep in certain cases a strong relation
to a Jacobian determinant -- so to say, an analogue of the relationship between principal curves and factors of the
classical {\em Jacobian curve} (see, e.g., Proposition~\ref{inversion_factor_is_determinant}).
Our main interest here in these inversion factors is the significant role they play as regards the generation
of some symbolic powers.

We succeed in going this far for general values $m\geq n\geq 3$.
To thrive deeper, we assume that either $m=n$ (the ``Cremona case'') or $m=n+1$ (the ``implicitization case'').
Our main drive is to tell the precise structure of the symbolic algebra $\mathcal{R}^{(I)}$ of $I$.
When $m=n$ our main results follows by drawing on some of the results of the earlier subsections
and collecting various pieces throughout the previous literature.
The main result says that $\mathcal{R}^{(I)}$ is generated in degrees $1$ and $n-1$, with only one fresh
generator in degree $n-1$ which may be taken to be the source inversion factor of the Cremona map
defined by the $n$-minors of $\mathcal{L}$.
Moreover, in characteristic zero, this generator coincides up to a scalar with the Jacobian determinant
of those same minors.

The case $m=n+1$ requires a full tour de force across the material and does not follow straightforwardly from the previously
stated results in the paper. First, the generation of $\mathcal{R}^{(I)}$ is more involved, occurring in degrees
$1$, $n-1$ and $n(n-1)-1$. This time around, showing that the source inversion factors constitute a minimal
set of fresh generators in degree $n-1$ is far from straightforward. Here we resorted to local duality
as applied to ${\rm H}^0_{(\XX)}(R/I^{n-1})\simeq I^{(n-1)}/I^{n-1}$ and to a subtle result on the $R$-dual
of the last nonfree syzygy module in the minimal free resolution of $R/I^{n-1}$. The argument here depends
strongly on the basic assumption that $I$ is the ideal of $n$-minors of a matrix whose entries are general
linear forms -- the result crumbles down for matrices with linear entries lacking this property.

This is the first step. In order to advance into proving the generation of the symbolic algebra we describe
a set of generators of its defining ideal, much in the spirit of \cite[Sections 5--8]{Wolmbook}, but  quite a bit more involved.
Making these generators explicit forced us to uncover a whole world of very tight relation between the various constructs
coming from the melange of symbolic power and birational theories.
A particular aspect that makes a case for this assertion is the long proof required to show that a certain variable is not
a zerodivisor modulo the ideal generated by the `expected' symbolic relations (proof of Theorem~\ref{symbolic_m=n+1}).
We have applied Gr\"obner basis theory via a case-by-case $S$-polynomial analysis in which the conclusions depend
strongly on the theoretical material developed before. Thus, it is not really the algorithm that matters, but rather
the use of the previous theory as a quality control.
Due to the amount of technical passages, we refer the reader to the appropriate places in the paper.

\section{Terminology}

%This section is devoted to collecting the basic elements of symbolic power theory
%and of birational theory as plied to our needs.

\subsection{Generalities on symbolic powers}

We will assume throughout that $R=k[X_1,\ldots,X_n]$ is a standard graded polynomial ring over an infinite field $k$.
Given an ideal $I\subset R$ and an integer $r\geq 1$, the $r$th {\em symbolic power} $I^{(r)}$ of $I$
is the contraction of $S^{-1}I^{r}$ under the natural homomorphism $R\rar S^{-1}R$
of fractions, where $S$ is the complementary set of the union of the associated primes of $R/I$.
In this work $I$ will be a codimension $2$ perfect ideal, hence $R/I$ is Cohen--Macaulay and so $I$ is a
pure (unmixed) ideal.
In this setup then $I^{(r)}$ is precisely the intersection of the primary components
of the ordinary power $I^{r}$ relative to the associated primes of $R/I$, i.e., the unmixed part of $I^{r}$.

A slightly different way to envisage symbolic powers is by noting that the $(I^{(r)}\cap I^{r-1})/I^{r}$
is the  $R/I$-torsion of the conormal module $I^{r-1}/I^{r}$ of order $r$.
Taking the direct sum over all $r\geq 0$ yields the $R/I$-torsion of the associated graded ring of $I$,
hence the non triviality of symbolic powers gives a measure of the torsion of the latter.
In particular, there is no nonzero torsion if and only if $I^{(r)}=I^r$ for every $r\geq 0$ -- in which
case one says that the ideal $I$ is {\em normally torsionfree}.
However, this information is most of the times pretty useless once it holds.
What matters for a substantial class of ideals -- codimension $2$ perfect ones included -- is
to guess some sort of asymptotic behavior for the equality of the two powers, more like an ``inf-asymptotic''
such behavior in the sense that one has equality throughout up to a certain exponent
order, thereafter comparison gets disorganized or even chaotic.

We observe that, like the ordinary powers, the symbolic powers constitute a decreasing multiplicative filtration,
so one can consider the corresponding {\em symbolic Rees algebra} $\mathcal{R}_{R}^{(I)}=\bigoplus_{r\geq 0}I^{(r)}t^{r}
\subset R[t]$.
However, unlike the ordinary Rees algebra, this algebra may not be finitely generated over $R$.
Alas, there are no definite effective ways to check when $\mathcal{R}_{R}^{(I)}$ is Noetherian.
The necessary and sufficient conditions of Huneke (\cite[Theorems 3.1 and  3.25]{Huneke}) obtained in dimension $3$
are not effective and neither is the necessary condition of Cowsik--Vasconcelos (\cite{Cowsik}, \cite[Proposition 3.5]{dual}).
Nevertheless, the latter becomes quite effective provided one has a good guess about what finitely generated subalgebra
looks like a strong candidate.
In a precise way, one has the following strategy.

First recall that, given an ideal $I\subset R$, where $R$ is a Noetherian domain with field of fractions $K$,
the {\em ideal transform} of $R$ relative to $I$ is the $R$-subalgebra $T_R(I):=R:_KI^{\infty}\subset K$.
We will draw on the following two fundamental facts:

\begin{itemize}
\item (\cite[Proposition 7.1.4]{Wolmbook}) If $C\subset T_R(I)$ is a finitely generated $R$-subalgebra such that
${\rm depth}_{IC}(C)\geq 2$ then $C=T_R(I)$.
\item (\cite[Proposition 7.2.6]{Wolmbook}) If $R$ moreover satisfies the condition $(S_2)$ of Serre then
$$\mathcal{R}^{(I)}_R\simeq T_{\mathcal{R}(I)}(J)\subset R[t]$$
as $R$-subalgebras of $R[t]$ for suitable choice of the ideal $J\subset R$.
\end{itemize}
Our idea of applying these principles is summarized in the following result,
of immediate verification:

\begin{Proposition}\label{idealtransform} Let $R=k[X_1,\ldots, X_n]$  denote a standard graded polynomial ring over an infinite field $k$,
with irrelevant maximal ideal  $(\XX):=(X_1,\ldots, X_n)$.
Let $I\subset R$ stand for a homogeneous ideal satisfying the following properties:
\begin{enumerate}
\item[{\rm(i)}] For every $r\geq 0$, the $R$-module $I^{(r)}/I^r$ is either zero or $(\XX)$-primary.
\item[{\rm(ii)}]  ${\rm depth}_{(\XX)C}(C)\geq 2$ for some finitely generated graded $R$-subalgebra $C\subset \mathcal{R}^{(I)}_R$
containing the Rees algebra $\mathcal{R}_R(I)$.
 \end{enumerate}
Then $C=\mathcal{R}^{(I)}_R$.
\end{Proposition}

We observe that the typical graded $R$-subalgebra $C\subset \mathcal{R}^{(I)}_R$
containing the Rees algebra $\mathcal{R}_R(I)$ as above has the form
$C=R[It,I^{(2)}t^2,\ldots,I^{(s)}t^{s}]\subset R[t]$, for suitable $s\geq 1$.
Although the non-vanishing of certain of the $R$-modules $I^{(r)}/I^r$ gives a measure of how far one has to go
(provided the symbolic Rees algebra is finitely generated), it is really the $R$-modules
$$\frac{I^{(r)}}{ \sum_{1\leq j\leq r-1}I^{(r-j)}\cdot I^{(j)}}$$
that count for the  search of {\em fresh} (or {\em genuine}) generators of the algebra.
Although this is a well-known simple observation, it often encrypts  some subtleties in a particular case.

\subsection{Birational maps and inversion factors}

Our reference for the basics in this part is \cite{bir2003}, which contains enough of the introductory
material in the form we use here (see also \cite{AHA} for a more general overview).

Let $k$ denote an arbitrary infinite field  -- further assumed to be algebraically closed in a geometric discussion.
A rational map $\mathfrak{G}:\pp^{n-1}\dasharrow \pp^{m-1}$ is defined by $m$ forms $\mathbf{g}=\{g_1,\ldots, g_m\}
\subset R:=k[\XX]=k[X_1,\ldots,X_n]$ of the same degree $d\geq 1$, not all null.
We naturally assume throughout that $n\geq 2$.
We often write $\mathfrak{G}=(g_1:\cdots :g_m)$ to underscore the projective setup
and assume that $\gcd\{g_1,\cdots ,g_m\}=1$ (in the geometric terminology, the linear system defining $\mathfrak{G}$ ``has no fixed part''),
in which case we call $d$ the {\em degree} of $\mathfrak{G}$.

Although the the definition of the rational map $\mathcal{G}$ depends on the linear system spanned by the defining coordinates,
its scheme theoretic indeterminacy locus is defined by the ideal of $R$ generated by the members of this system.
For convenience, this ideal will slightly improperly be referred to as the {\em base ideal} of $\mathfrak{G}$.

The image of $\mathfrak{G}$ is the projective subvariety $W\subset \pp^{m-1}$ whose homogeneous
coordinate ring is the $k$-subalgebra $k[\mathbf{g}]\subset R$ after degree renormalization.
Write $k[\mathbf{g}]\simeq k[\YY]/I(W)$, where $I(W)\subset k[\YY]=k[Y_1,\ldots,Y_m]$ is the homogeneous defining ideal
of the image in the embedding $W\subset \pp^{m-1}$.

We say that $\mathfrak{G}$ is {\em birational onto the image} if there is a rational map
 $\pp^{m-1}\dasharrow \pp^{n-1}$ with defining coordinate forms $\mathbf{f}=\{f_1,\ldots, f_n\}
\subset k[\YY]$ (not simultaneously vanishing modulo $I(W)$)  satisfying the
relations
\begin{equation}\label{birational_rule}\nonumber
({f}_1(\mathbf{g}):\cdots :{f}_n(\mathbf{g}))=(X_1:\cdots :X_n), \;
({g}_1(\mathbf{f}):\cdots :{g}_m(\mathbf{f}))\equiv (Y_1:\cdots :Y_m)\pmod {I(W)}
\end{equation}
Let $K$ denote the field of fractions of $k[\mathbf{g}]$.
The coordinates $\{f_1,\cdots ,f_n\}$ defining the ``inverse'' map are not uniquely defined;
any other set $\{f'_1,\cdots ,f'_n\}$  inducing the same
element of the projective space $\pp^{n-1}_{K}=\pp^n_k\otimes_k {{\rm Spec}(K)}$ will do as well -- both tuples are called
{\em representatives} of the rational map. Furthermore, one can choose a finite minimal set
$\mathbf{f_1},\ldots, \mathbf{f_s}$ of these representatives
such that any other representative belongs to the $k[\YY]/I(W)$-submodule generated by
$\mathbf{f_1},\ldots, \mathbf{f_s}$. More exactly, any such a minimal representative is the transpose
of a minimal generator of the syzygy module of the so-named weak Jacobian dual matrix
(for complete details see \cite{bir2003}, particularly Proposition 1.1 and \cite[Section 2]{AHA}).
Such a set will be referred to in the sequel as {\em a complete set of minimal representatives} of the
inverse map.

Having information about the inverse map -- e.g., about its degree -- will be quite relevant in the sequel.
For instance, the first of the above structural congruences
\begin{equation}\label{congruence}
(f_1(g_1,\ldots,g_m),\ldots, f_n(g_1,\ldots,g_m))\equiv (X_1,\ldots,X_n)
\end{equation}
involving the inverse map, in terms of a given representative lifted to $k[\yy]$,
yields a uniquely defined
form $D\in R$ up to a nonzero scalar in $k$, such that $f_i(g_1,\ldots,g_m)=X_iD$, for every $i=1,\ldots,n$.
This is merely a consequence of factoriality. Indeed, the congruence means that there are forms $D,D'$ such
that $D' f_i(\mathbf{g})= Dx_i$ for every
$1\leq i\leq n$.  Now, a prime factor of $D'$ that does not divide $D$ would have to divide all $n\geq 2$
variables, which is only possible if $D'\in k$. Otherwise, necessarily $D'$ divides $D$; in any case we are through.

We call $D$ the {\em source inversion factor} of $\mathfrak{G}$ associated to the given representative.
There is a counterpart {\em target inversion factor}, defined in a parallel way by exchanging the roles of
$\ff$ and $\gg$ in (\ref{congruence}).
However, for that to be seen one has to be slightly more careful, as we now explain.
For convenience we state the result in the form of a lemma.

\begin{Lemma}\label{uniqueness_of_target_factor}
Let $\mathfrak{G}:\pp^{n-1}\dasharrow \pp^{m-1}$ be a rational map defined by forms $\mathbf{g}=\{g_1,\ldots, g_m\}
\subset R:=k[\XX]=k[X_1,\ldots,X_n]$ of the same degree.
Suppose that $\mathfrak{G}$ is birational onto a non-degenerate subvariety $W\subset \pp^{m-1}$.
Then the inverse map admits a representative by forms $\mathbf{f}=\{f_1,\ldots, f_n\}
\subset k[\YY]$ of the same degree such that there exists a uniquely defined form $E\subset k[\YY]$ modulo $I(W)$ satisfying the congruences
$$(g_1(\mathbf{f}),\cdots ,g_m(\mathbf{f}))\equiv E\,\cdot (Y_1,\cdots ,Y_m)\pmod {I(W)}.$$
\end{Lemma}
\demo The proof is strongly based on the results and proof of \cite[Theorem 2.18]{AHA}.
We follow verbatim the line of  argument of the proof of \cite[Theorem 2.18, Supplement, (ii)]{AHA}.
First,  a representative of the inverse map $\mathbf{f}=\{f_1,\ldots, f_n\}$ can be taken in which $f_i$
denotes a signed $(n-1)$th order minor of an $(n-1)\times n$ submatrix of the weak Jacobian dual matrix of $\mathfrak{G}$.
Next, drawing upon the so-called Koszul--Hilbert Lemma (\cite[Proposition 2.1]{AHA}) one derives a set of simpler congruences:
\begin{equation}\label{from_KH}
g_m(\mathbf{f})\, Y_j =  g_j(\mathbf{f})\, Y_m  \pmod {I(W)},\; 1\leq j\leq m.
\end{equation}
Since $W\subset \pp^{m-1}$ is non-degenerate, the ideal $(Y_m, I(W))\subset k[\YY]$
is prime.
We read the above congruences as $g_m(\mathbf{f})\, Y_j\in (Y_m, I(W))$, for $1\leq j\leq m$.
If $g_m(\mathbf{f})\not\in (Y_n, I(W))$ then $Y_j\in (Y_m, I(W))$ for $1\leq j\leq m$.
This gives $(Y_1,\ldots, Y_m)\subset (Y_n, I(W))$ which implies that $I(W)$ has codimension $m-1$, hence
$\dim W=0$. But this is impossible since the map is birational and $n\geq 2$.
Therefore, one must conclude that $g_m(\mathbf{f})\in (Y_m, I(W))$, and hence $Y_m$ is cancellable in (\ref{from_KH}).
\qed

\smallskip

We remark that $E$ depends on the choice of the forms $\{f_1,\ldots, f_n\}$ representing the
inverse map, which is by no means uniquely given (see \cite{AHA} for the details of this nature).

\medskip

A fundamental property of the inversion factor of a Cremona map in characteristic zero does not 
seem to have been observed before in the following generality and explicitness.
We give a neat algebraic proof.

\begin{Proposition}\label{jac_vs_factor}{\rm (char$(k)=0)$}
Let $\mathfrak{G}$ denote a Cremona map of $\pp^{n-1}$ defined by forms $\gg:\{g_1,\ldots,g_n\}$ in $R$
without fixed part and let $\Theta(\gg)$ denote the Jacobian matrix of $\gg$.
Then $\det (\Theta(\gg))$ divides a power of the source inversion factor $G$ of $\mathfrak{G}$.
In particular, if $\det (\Theta(\gg))$ is reduced then it divides $G$.
\end{Proposition}
\demo
Let $\ff:\{f_1,\ldots,f_n\}$ define the inverse map.
Applying the chain rule of derivatives to the structural equation
${\bf f}({\bf g})=G\cdot(\XX)$, it obtains
\begin{equation}\label{derivatives}
\Theta({\bf f})({\bf g})\cdot\Theta({\bf g})=G\cdot \mathcal{I}+(\XX)^t\cdot {\rm Grad}(G)
\end{equation}
where $\mathcal{I}$ is the identity matrix and ${\rm Grad}(G)=(\partial G/\partial X_1 \ldots \partial G/\partial X_n)$.
Note that the right side of (\ref{derivatives}) is the result of evaluating $\lambda\mapsto G$ in the characteristic matrix
$\lambda\mathcal{I}-\mathcal{A}$, where $\mathcal{A}=-(\XX)^t\cdot {\rm Grad}(G)$.

Recall that, quite generally the coefficients of the characteristic polynomial $p(\lambda)=\lambda^n+a_1\lambda^{n-1}+\ldots+a_{n-1}\lambda+a_n$ of $\mathcal{A}$ can be recursively computed as:
 \begin{eqnarray}\label{leverrier}
-a_1&=&s_1\nonumber\\
-ra_r&=&s_r+\sum_{i=1}^{r-1}s_{i}a_{r-i}
 \end{eqnarray}
where $s_r$ is the trace of the  matrix $\mathcal{A}^r$, for $1\leq r\leq n.$

Now, one has
 \begin{equation}
  \mbox{Trace}((\XX)^t\cdot {\rm Grad}(G))=\sum_{i=1}^nX_i\frac{\partial G}{\partial X_i}=d\,G,
 \end{equation}
 where $d=\deg(G)$.
 On the other hand, a calculation yields the equality $\mathcal{A}^2=(d\, G) \mathcal{A}$.
 By an immediate recursion it follows that
 \begin{equation}
 ((\XX)^t\cdot {\rm Grad}(G))^r=(d^{r-1}G^{r-1})(\XX)^t\cdot {\rm Grad}(G)
\end{equation}
From this, applying (\ref{leverrier}) recursively yields $a_2=a_3=\cdots =a_n=0$, hence $p(t)=t^n+(d\,G)t^{n-1}$.
Evaluating $\lambda\mapsto G$ yields $p(G)=G^n+(d\,G)G^{n-1}=(d+1)G^n$.
But this is $(\det G\cdot \mathcal{I}+(\XX)^t\cdot {\rm Grad}(G))$, hence
(\ref{derivatives}) gives
\begin{equation}\label{fator}
\det(\Theta({\bf f})({\bf g}))\cdot\det(\Theta({\bf g}))=(d+1)G^{n}.
\end{equation}
Therefore, $\det(\Theta({\bf g}))$ indeed divides $G^n$.
\qed

\medskip

From the other end, the basic result relating inversion factors to symbolic powers is
the following proposition showing how the former become genuine generators of the latter.
We are not aware of this result having been explicitly pointed out
in the previous literature.

\begin{Proposition}\label{inversionfactor_is_symbolic}
Let $I\subset R=k[\xx]$ denote the base ideal of a birational map $\mathfrak{F}:\pp^{n-1}\dasharrow \pp^{m-1}$ onto
the image, satisfying the canonical restrictions. Let $D\subset R$ denote the source inversion factor relative
to a given minimal representative of the inverse map.
Suppose that $I$ is a radical ideal such that $I^{(\ell)}=I^{\ell}$ for every $\ell\leq d'-1$,
but $I^{(d')}\neq I^{d'}$, where $d'$ is the degree of the coordinates of the representative.
Then
\begin{enumerate}
\item[{\rm (a)}] $D\in I^{(d')}\setminus \sum_{1\leq j\leq d'-1}I^{(d'-j)}\cdot I^{(j)}=I^{d'}$.
\item[{\rm (b)}] Moreover, if $I^{(d')}$ is generated in standard degree $\geq dd'-1$, where $d$ is the degree of
the coordinates of $\mathfrak{F}$, then $D$ is a homogeneous minimal generator of the symbolic Rees algebra.
 \end{enumerate}
\end{Proposition}
\demo
The characteristic property of $D$ is the congruence (\ref{congruence}).
In particular, $D\in I^{d'}: (\xx)$.
We may assume that $\mathfrak{F}$ is not the identity map of $\pp^{n-1}$.
Since $I$ is radical, it has codimension at most $\dim R-1$; hence, there is a form $h\in (\xx)\setminus P$,
for every minimal prime $P$ of $R/I$,  such that $hD\in I^{d'}$.
This means that $D\in I^{(d')}$.

The rest of (a) follows straightforwardly
under the present hypothesis.

Part (b) is clear since $\deg(D)=dd'-1$.
\qed

\section{Ideals of general linear forms}

\subsection{Arithmetic and homological properties}

Let $k$ stand for an infinite field  and let $R=k[X_1,\ldots, X_n]$ denote a
standard graded polynomial ring over $k$.
We will often require that char$(k)=0$, but some of the results will be valid in any characteristic.
Our basic object is an $m\times (m-1)$ matrix of general linear forms over $R$.
For the sake of subsequent ideal theoretic properties, we introduce our main notion in a
more formal way, by emphasizing its nature as a specialization out of a {\em generic} $m\times (m-1)$ matrix
over a larger polynomial ring.

%Let $\{F_1,\ldots,F_s\}$ be a set of generic forms over $k$ of
%respective degrees $d_1,\ldots, d_s$.
%More precisely, we have $F_j=\sum_{|a|=d_j} T_{j,a}\YY^a,$
%$where $\TT:=\{T_{j,a},\,1\leq j\leq s,\, |a|=d_j\}$ are mutually independent
%indeterminates over $k$ and $\YY=\{Y_1,\ldots,Y_p\}$ is a set of $p$  variables
%over $k[\TT]$.
%If $s<p$ then $\{F_1,\ldots,F_s\}$ is a regular sequence (see, e.g., \cite[Proposition 4.3]{Jou}).
%Now map the variables  $T_{j,a}\mapsto \alpha_{j,a}\in k$, where the $\alpha_{j,a}$'s are random elements.
%The resulting set of forms $f_1,\ldots,f_s\in k[\YY]$ will still be a regular sequence
%of forms, but with random coefficients.
%Such a set will be called a {\em general regular sequence} of forms.
%We emphasize that this means that the entire set of coefficients of the $f_i$'s is general - not merely
%within an individual form.

\begin{Definition}\label{DEF}\rm
Let $m, n$ be given positive integers such that $m(m-1)\geq n$.
Let $\mathcal{Z}=(Z_{i,j})$ denote an $m\times (m-1)$ generic matrix over $k$.
Write $S:=k[\mathbf{Z}]=k[Z_{i,j}\,|\, 1\leq i\leq m, 1\leq j\leq m-1]$ for the corresponding polynomial ring over $k$.
Setting $s:=m(m-1)-n$, let ${\mathbf{L}}:=(L_1,\ldots,L_s)\subset S$ denote a regular sequence of general $1$-forms
(i.e., the entire set of coefficients is random) -- such regular sequences abound by taking a regular
sequences of $s$ forms of degree $1$ with indeterminate coefficients (see, e.g., \cite[Proposition 4.3]{Jou}))
and then randomly evaluating these indeterminates to elements of $k$.
Denote $R:=S/(\mathbf{L})$ and let $\mathcal{L}$ stand for the $m\times (m-1)$ matrix over $R$ whose entries are the images of
the variables $Z_{ij}$ under the surjection of $k$-algebras $S\surjects R$.
For lack of better terminology, we call $\mathcal{L}$  a {\em general linear} matrix
of size $m\times (m-1)$.
\end{Definition}
Clearly then $R$ is isomorphic to the ring of polynomials over $k$ on $n$ variables.
For the sake of subsequent development, we make this setup more explicit in the following way.
Namely, consider the revlex monomial order of $S$ with the following ordering of the variables:
$Z_{1,1}>\cdots >Z_{1,m-1}> Z_{2,1}\cdots >Z_{m, m-1}$.
It is a simple inductive argument on $s$ using elementary operations on the generators to see that the
first $s=m(m-1)-n$ variables in this order
generate an ideal contained in the initial ideal of $\mathbf{L}$.
But since the former is a prime ideal of codimension $s$ it must coincide with the latter.
Likewise, the inductive procedure used also gives that
 $\mathbf{L}$ is generated by  $1$-forms $L_{i,j}:=Z_{i,j}-\lambda_{i,j}$, where $\{{i,j}\}$ runs
through the first $m(m-1)-n$ entry indices and the corresponding $\lambda_{i,j}$
is an $1$-form depending only on the last $n$ variables.
One notes that the entire set of coefficients of the set $\{\lambda_{i,j}\,|\, i,j\}$ of $1$-forms
is a result of simple operations on the original coefficients of the set $\{L_1,\ldots,L_s\}$, hence is itself
general, if not strictly random.
Then $R$ is isomorphic to the ring of polynomials over $k$ on these variables.
We now rename these variables to $X_1,\ldots,X_n$, and henceforth write $R=k[X_1,\ldots,X_n]$.
With this proviso, the matrix $\mathcal{L}$ is an $m\times (m-1)$ matrix over $R$ where the first $m(m-1)-n$
entries (in the entry ordering as above) are the forms $\ell_{ij}$ obtained by evaluating the forms
$\lambda_{i,j}$  on $X_1,\ldots,X_n$ and the last $n$ entries are the variables $X_1,\ldots,X_n$ themselves.

One can harmlessly trade the last $n$ entries for $n$ additional random linear forms in $X_1,\ldots,X_n$.
We emphasize once more that the entire set of the forms $\ell_{i,j}$ appearing as entries of
the matrix is general, i.e., the totality of all coefficients is random.

Following common usage, one denotes by $I_t(\Psi)\subset R$ the ideal generated by the $t\times t$ minors of a matrix $\Psi$.
In the present setup, one has $I_t(\mathcal{L})=(I_t(\mathcal{Z}),\mathbf{L})/(\mathbf{L})$.

Our first basic result is about the codimension of these ideals of minors.

\begin{Theorem}\label{specialization}
Let $\mathcal{L}$ denote an $m\times (m-1)$
general linear matrix over $R=k[X_1,\ldots, X_n]$ as above, with $n\geq 3$ and $m\geq 2$.
Then, for every $1\leq t\leq m-1$, one has
$$\hht(I_t(\mathcal{L}))= \min\{n, (m-t+1)(m-t)\}.$$
\end{Theorem}
\demo
Because $S/(I_t(\mathcal{Z}),\mathbf{L})\simeq R/I_t(\mathcal{L})$, the assertion is equivalent
to the following one
\begin{equation}\label{max_dim}
\dim S/(I_t(\mathcal{Z}),\mathbf{L})=\max\{0, n- (m-t+1)(m-t)\}.
\end{equation}
 More precisely, we will now show that, for $1\leq t\leq m-1$ and $1\leq r\leq s=m(m-1)-n$:
\[
\dim S/(I_t(\mathcal{Z}),L_1,\ldots,L_{r})=\left\{
\begin{array}{ll}
0 & \mbox{if $D < r$}\\
D-r & \mbox{if $D\geq r$}
\end{array}
\right.
\]
where $D:=\dim S/I_t(\mathcal{Z})=m(m-1)- (m-t+1)(m-t)$.

We proceed by induction on $r$.
Obviously, $D>0$ if and only if $t\geq 2$.
Now, for every $t$ in this range, clearly $L_1$ is a non-zero-divisor on $S/I_t(\mathcal{Z})$ since $L_1$ is a linear form
and all the associated primes are contained in the one single prime $I_2(\mathcal{Z})$ generated in degree $2$.
Therefore, one has

\[
\dim S/(I_t(\mathcal{Z}),L_1)=\left\{
\begin{array}{ll}
0 & \mbox{if $D=0$}\\
D-1 & \mbox{if $D\geq 1$}
\end{array}
\right.
\]
Let now $m(m-1)-n\geq j\geq 2$.
By the inductive hypothesis, one has
\[
\dim S/(I_t(\mathcal{Z}),L_1,\ldots,L_{j-1})=\left\{
\begin{array}{ll}
0 & \mbox{if $D< j-1$}\\
D-(j-1) & \mbox{if $D\geq j-1$}
\end{array}
\right.
\]
Now consider the set $\Ass(S/I_t(\mathcal{Z}),L_1,\ldots,L_{j-1})$, with $t$ in the range
for which $D\geq j-1$. This is a finite set of primes.
Let $J_{[1]}$ denote the part of degree $1$ of a homogeneous ideal $J$ in $R$.
Since $L_j$ is a general form and randomly chosen with respect to the forms $L_1,\ldots,L_{j-1}$, we have
$L_j\notin P_{[1]}$ for every prime $P\in \bigcup_{t}\Ass(S/I_t(\mathcal{Z}),L_1,\ldots,L_{j-1})$.

Then the dimension again drops by $1$, i.e., we get
\[
\dim S/(I_t(\mathcal{Z}),L_1,\ldots,L_{j})=\left\{
\begin{array}{ll}
0 & \mbox{if $D< j$}\\
D-j & \mbox{if $D\geq j$}
\end{array}
\right.
\]
Applying with $j=m(m-1)-n$ yields
\[
\dim S/(I_t(\mathcal{Z}),\mathbf{L})=\left\{
\begin{array}{ll}
0 & \mbox{if $D < m(m-1)-n$}\\
D-m(m-1)+n & \mbox{if $D\geq m(m-1)-n$}
\end{array}
\right.
\]
Substituting for the value of $D$ yields the required result.
\qed

\medskip

Next is the first basic structural result.

\begin{Proposition}\label{general_cod2}
Let $\mathcal{L}$ denote an $m\times (m-1)$
general linear matrix over $R=k[X_1,\ldots, X_n]$, with $n\geq 3$ and $m\geq 2$.
Set $I=I_{m-1}(\mathcal{L})\subset R$.
Then:
\begin{itemize}
\item[{\rm (i)}] $I$ has codimension $2$ and $I_{m-2}(\mathcal{L})\subset R$  has codimension $\min\{6,\, n\}$.
\item[{\rm (ii)}] {\rm (char}$(k)=0${\rm )} $R/I$  satisfies the condition $(R_r)$ of Serre, with
$r=\min\{3,\, n-2-1\}\,${\rm ;} in particular, if $n\geq 4$ then
$R/I$ is normal and $I$ is a prime ideal.
\item[{\rm (iii)}] $I$ is of linear type  if and only if
$m\leq n$.
\item[{\rm (iv)}] {\rm (char}$(k)=0${\rm )} $I$ is normally torsionfree if and only if $m<n$.
\end{itemize}
\end{Proposition}
\demo
(i) This follows from Theorem~\ref{specialization}.

\smallskip

(ii)
Since $\mathcal{Z}$ is a generic matrix,  the Jacobian ideal of $S/I_{m-1}(\mathcal{Z})$ is
$I_{m-2}(\mathcal{Z})/I_{m-1}(\mathcal{Z})$.
Applying Bertini's theorem (\cite{FlenVog}) gives that the singular scheme of the scheme-theoretic
general hyperplane section $S/(I_{m-1}(\mathcal{Z}),L_1)$ is the scheme associated to $S/(I_{m-2}(\mathcal{Z}),L_1)$.
Inducting on the number $m(m-1)-n$ of general hyperplane sections yields that the singular scheme of the scheme-theoretic
linear section $S/(I_{m-1}(\mathcal{Z}),\mathbf{L})\simeq R/I$ is the scheme associated to
$S/(I_{m-2}(\mathcal{Z}),\mathbf{L})\simeq R/I_{m-2}(\mathcal{L})$.
By Theorem~\ref{specialization}, the latter has codimension
 at least $\min\{6,\, n\}$ on $R$.
Since $I=I_{m-1}(\mathcal{L})$ has codimension $2$, $R/I$ satisfies  the Serre condition $(R_{\min\{3,n-2-1\}})$.
Thus,  if $n\geq 4$ then $R/I$ satisfies $(R_1)$.
At the other end, $R/I$ is Cohen--Macaulay.
It follows that, for $n\geq 4$, $R/I$ is normal and, since $I$ is homogeneous, $R/I$ must be a domain.
(If $n=3$ then $I$ is still a radical ideal.)

\smallskip

(iii) Let us apply the result of Theorem~\ref{specialization} in this case.
We claim that
$$\min\{n, (m-t+1)(m-t)\}\geq m-t+1,\; {\rm for}\; 1\leq t\leq m-1.$$
This is obvious if the minimum is attained by $(m-t+1)(m-t)$; if the minimum is $n$ instead then
$m\leq n$ certainly implies $m-t+1\leq n$.

This shows that $I$ satisfies the property $(F_1)$, hence it is an ideal of linear type in
this case (see \cite{HSV1}).
The converse is evident since the linear type property implies the inequality $\mu(I)\leq \dim R$.

\smallskip

(iv) Suppose first that $m<n$.
By part (iii), $I$ is of linear type. Since $I$ is strongly Cohen--Macaulay (\cite[Theorem 2.1(a)]{AVRAMOVHERZOG}) then
the Rees algebra of $I$ is Cohen--Macaulay (\cite[Theorem 9.1]{HSV1}), and hence so is the associated graded ring of $I$.
On the other hand, we may assume that $n\geq 4$ given that for $n=3$ the ideal $I$ is generated by a regular sequence of two
elements.
Therefore, by part (ii), the ideal $I$ is prime. By \cite[Proposition 3.2 (1)]{EH}, the assertion is equivalent to having
$$ \ell_P(I)\leq \max\{\hht P-1, \hht I\},$$
for every prime ideal $P\supset I$. Since $I$ is homogeneous, it suffices to take $P$ homogeneous.
We may assume that $\hht P\geq 3$ since $I$ is a height $2$ prime. Therefore, we have to show that $\ell_P(I)\leq \hht P-1$.
If $P=(\XX)$ the result is clear since $\ell_{(\XX)}\leq \mu(I)=m\leq n-1=\hht (\XX)-1$.
Therefore, we may assume that $P\subsetneq (\XX)$, hence $\hht P\leq n-1$.

Set $t_0:=\max\{1\leq s\leq m-2\,|\, I_s(\mathcal{L})\not\subset P\}$.
Therefore, $I_{t_0+1}(\mathcal{L})\subset P$, hence $\hht I_{t_0+1}(\mathcal{L})\leq \hht P\leq n-1$.
By Theorem~\ref{specialization} one must have $\hht I_{t_0+1}(\mathcal{L})=(m-t_0)(m-t_0-1)$.
Pick a $t_0$-minor  $\Delta$ of $\mathcal{L}$ not contained in $P$, so that, in particular, $R_P$ is
a localization of the ring of fractions $R_{\Delta}=R[\Delta^{-1}]\subset k(\XX)$.
By a standard row-column elementary operation procedure, there is an $(m-t_0)\times(m-t_0-1)$ matrix $\widetilde{\mathcal{L}}$ over
$R_P$ such that
$$I_P=I_{m-1-t_0}(\widetilde{\mathcal{L}}).$$
Assume first that $t_0\leq m-3$.
Then $(m-t_0)\leq (m-t_0)(m-t_0-1)-1=\hht I_{t_0+1}(\mathcal{L})-1\leq\hht P -1$.
Therefore
\begin{equation*}
\ell_P(I)=\ell (I_{m-1-t_0}(\widetilde{\mathcal{L}}))\leq\min\{\mu(I_{m-1-t_0}(\widetilde{\mathcal{L}})), \hht P\}
=\min\{m-t_0, \hht P\}\leq\hht P-1.
\end{equation*}
If $t_0=m-2$,  one gets $\ell_P(I)=\min\{2, \hht P\}=2\leq \hht P-1$ since it has been assumed that $\hht P\geq 3$.

Therefore, $I$ is normally torsionfree.
The converse will follow from Theorem~\ref{general_symbolic_cod2}.
\qed

\medskip

The proof of the main theorem stated further down  will draw on several results of independent interest.

Recall the following notation: for a given integer $s\geq 1$, one says that the ideal $I\subset R$ satisfies
condition $(G_s)$ if $\mu(I_P)\leq \hht P$, for every prime ideal $P$ such that $\hht P\leq s-1$.

\begin{Proposition}\label{alma_mater}
Let $\mathcal{L}$ denote an $m\times (m-1)$ general linear matrix over $R=k[X_1,\ldots, X_n]$, with $m \geq 2$ and $n\geq 3$.
Set $I:=I_{m-1}(\phi).$ Then
\begin{enumerate}
\item[{\rm (i)}] $I$ satisfies condition $(G_n)$.
\item[{\rm (ii)}] {\rm (char}$(k)=0${\rm )} Given an integer $r\geq 0$ such that $I^{(r)}/I^r\neq \{0\}$ then $I^{(r)}/I^r$ is an
$(\XX)$-primary $R$-module {\rm (}in other words, $I^{(r)}$ is the saturation of $I^r${\rm )}.
\end{enumerate}
\end{Proposition}
\demo
(i) Let $P\subset R$ be a prime of height $\leq n-1$.
Set
$$t_{\infty}:=\min\{1\leq t\leq m-1\,|\, I_t(\mathcal{L}\subset P)\}.$$
Then $\hht I_{t_{\infty}}(\mathcal{L})\leq n-1$, hence $\hht I_{t_{\infty}}(\mathcal{L})=(m-t_{\infty}+1)(m-t_{\infty})$
by Proposition~\ref{general_cod2} (i).
Inverting a $(t_{\infty}-1)$-minor of $\mathcal{L}$ in $R_P$ we  get
$I_P=I_{m-t_{\infty}}(\widetilde{\mathcal{L}})$ for a suitable $(m-t_{\infty}+1)(m-t_{\infty})$ matrix
$\widetilde{\mathcal{L}}$ over $R_P$.
Collecting the information yields
\begin{eqnarray*}
\mu(I_P)&=& \mu(I_{m-t_{\infty}}(\widetilde{\mathcal{L}}))=m-t_{\infty}+1\leq (m-t_{\infty}+1)(m-t_{\infty})\\
&=& \hht I_{t_{\infty}}(\mathcal{L})\leq \hht(P).
\end{eqnarray*}

\smallskip

(ii) Fixing an $r\geq 0$, suppose that $I^{(r)}/I^r\neq \{0\}$.
By Proposition~\ref{general_cod2} (iv), we have $m\geq n$.
The assertion is equivalent to saying that a power of $(\XX)$ annihilates $I^{(r)}/I^r$
i.e., that ${I^{(r)}}_P={I^r}_P$ for every prime $P\neq (\XX)$.
Letting $r\geq 0$ run, this is in turn equivalent to claiming that the associated graded ring gr$_I(R)$ is torsionfree
over $R/I$ locally on the punctured spectrum Spec$(R)\setminus (\XX)$.

Thus, let $P\neq (\XX)$ be a prime containing $I$.
Then the condition $(G_n)$ of part (i) implies that $I_P$ satisfies the condition $(F_1)$ (same as $(G_{\infty})$)
as an ideal of $R_P$.
As in the proof of Proposition~\ref{general_cod2} (iv), we know that the associated graded ring gr$_{I_P}(R_P)$
is Cohen--Macaulay. Therefore, by the same token and since $\hht I=2$, one has to show the local estimates
$$\ell_Q(I)=\ell_{Q_P}(I_P)\leq \hht (Q_P)-1= \hht Q-1,$$
for every prime $Q\subset P$.

Fixing such a prime $Q$, set $t_0:=\max\{1\leq s\leq m-2\,|\, I_s(\mathcal{L})\not\subset Q\}$.
The argument is now the same as the one in the proof of Proposition~\ref{general_cod2} (iv).
\qed

\begin{Corollary}\label{F_1}
Let $\mathcal{L}$ denote an $m\times (m-1)$ general linear matrix over $R=k[X_1,\ldots, X_n]$,  with $m \geq 2$ and $n\geq 3$.
Set $I:=I_{m-1}(\phi).$ Then
 $\mathcal{S}_r(I)\simeq I^{r}$ in the range $1\leq r\leq n-1$.
\end{Corollary}
\demo
This follows from Proposition~\ref{alma_mater} (i) as applied through the result of \cite[Theorem 5.1]{TCHERNEV}.
\qed

\medskip

For an integer in the range $1\leq r\leq n-3$, recall the $r$th approximation complex associated to the ideal $I$
(see \cite[Section 3]{Wolmbook}):
\begin{equation}\label{M-complex}
\mathcal{M}_r:0 \rar H_r\rar H_{r-1}\otimes S_1\rar \cdots \rar H_1\otimes S_{r-1}\rar S_r.
\end{equation}

Here $H_i$ stands for the $i$th Koszul homology module on the generators of $I$ and $S_i$ denotes the $i$th
homogeneous part of the polynomial ring $S:=R/I[Y_1\ldots,Y_m].$
One has $H_0(\mathcal{M}_r)\simeq \mathcal{S}_r(I/I^2)$.

\begin{Proposition}\label{aciclicidade}
The approximation complex $\mathcal{M}_r$ is acyclic in the range  $1\leq r\leq n-3$.
\end{Proposition}
\demo
We show that the complex is acyclic locally everywhere.
Suppose first that $P\neq (\XX)$ is a non-irrelevant prime.
In this case, using Corollary~\ref{F_1}, the result is contained in \cite[Theorem 5.1]{HSV1}.

Thus, we can assume that $P= (\XX)$ and that $\mathcal{M}_r$ is acyclic locally at any prime
properly contained in $(\XX)$.
We show acyclicity stepwise from the left.
Thus, suppose the partial complex
$$\begin{array}{ccccccccc}
0 &\rar & H_r &\rar &\cdots &\rar & H_{k+2}\otimes_{r-k-2} & \rar & H_{k+1}\otimes S_{r-k-1}\\
&&&&&&&& \searrow \\
&&&&&&&& \;\;\quad B_k\\
 &&&&&&&& \;\;\quad\quad\quad \searrow\\
 &&&&&&&& \;\;\quad \quad\quad\quad\quad 0
\end{array}
$$
is exact.
Since $I$ is a strongly Cohen--Macaulay ideal (\cite[Theorem 2.1(a)]{AVRAMOVHERZOG}),  one has $\depth (H_i)=n-2$
for every $1\leq i\leq r$.
Chasing depths from left to right, one gets $\depth (B_k)\geq n-(r-k+1)=(n-r) +k-1\geq 3+k-1=k+2\geq 2$.

Now, letting $Z_k\subset H_{k}\otimes S_{r-k}$ denote the subsequent module of cycles, write $D_k:=Z_k/B_k$.
Suppose $D_k\neq 0$ and take $Q\in \Ass(D_k)$.
Since the entire complex is acyclic off $(\XX)$, we must have $Q=(\XX)$.
Applying Hom$_R(R/(\XX), -)$ yields the exact sequence
$$0=\Hom_R(R/(\XX),Z_k)\rar\Hom_R(R/(\XX),D_k)\rar\Ext {1}{R}{R/(\XX)}{B_k}.$$
The rightmost term of this sequence vanishes as well since the depth of $B_k$ is at least $2$, hence also does the middle term;
this is absurd since $(\XX)$ is an associated prime of $D_k$.
Therefore, we conclude that $D_k=0$.
\qed

\medskip

Denote by pd$_R(M)$ the projective dimension of a finitely generated $R$-module $M$.

\begin{Corollary}\label{dproj}
${\rm pd}_R(I^{r}/I^{r+1})\leq r+2$ for $1\leq r\leq n-3$. In particular, $(\XX)\notin\Ass(I^{r}/I^{r+1})$.
\end{Corollary}
\demo Since (\ref{M-complex}) is acyclic by Proposition~\ref{aciclicidade}, depth chasing all the way to the right
yields $\depth \mathcal{S}_r(I/I^2)\geq n-(r+2)$. Therefore, pd$_R(\mathcal{S}_r(I/I^2))\leq r+2$.
But $\mathcal{S}_r(I/I^2)\simeq I^{r}/I^{r+1}$ by Corollary~\ref{F_1}.
\qed

\begin{Theorem}\label{associated_primes} {\rm (char}$(k)=0${\rm )}
$\Ass(R/I^{r})=\Ass(R/I)$ for $1\leq r\leq n-2$.
\end{Theorem}
\demo
By Proposition~\ref{alma_mater} (ii), $\Ass(R/I^{r})\subset\Ass(R/I)\cup \{(\XX)\}$ - note that
the assumption that $I$ is prime holds for $n\geq 4$; for $n=3$, $I$ is still radical, hence the statement
is obvious directly.

Proceed by induction on $r$.
It is clear for $r=1$ since $I$ is a radical unmixed ideal for $n\geq 3$ and $\hht (I)=2<n$.

Supposing $(\XX)\in \Ass(R/I^{r})$, the exact sequence $0\rar I^{r-1}/I^{r}\rar R/I^{r}\rar R/I^{r-1}\rar 0$ and the inductive
hypothesis  force us to conclude that $(\XX)\in \Ass(I^{r-1}/I^{r})$. But since $r+1\leq n-1$, the latter is forbidden by
 Corollary~\ref{dproj}.
 \qed

\subsection{The role of the inverse factor}

An ideal $I\subset R$ generated by $m$ forms of the same degree  is {\sc of fiber type} if the
bihomogeneous defining ideal $\mathcal{J}\subset R[\YY]=R[Y_1,\ldots,Y_m]$ of the Rees algebra $\mathcal{R}(I)$
is generated by its $\YY$-linear forms and the defining equations of the special fiber
$\mathcal{R}(I)/(\XX)\mathcal{R}(I)$.

\begin{Proposition}\label{m_greater_than_n_is _birational}
Let $\mathcal{L}$ denote an $m\times (m-1)$
general linear matrix over $R=k[X_1,\ldots, X_n]$, with $m\geq n\geq 3$.
Setting $I:=I_{m-1}(\mathcal{L})\subset R$, one has:
\begin{itemize}
\item[{\rm (a)}] The rational map $\mathfrak{G}:\pp^{n-1}\dasharrow \pp^{m-1}$ defined by the $(m-1)$-minors
is birational onto its image.
\item[{\rm (b)}] $I$ is an ideal of fiber type and the Rees algebra $\mathcal{R}(I)$ is a Cohen--Macaulay domain.
\item[{\rm (c)}] The map $\mathfrak{G}$ admits ${{m-1}\choose {n-1}}$ source inversion factors, each associated
to a minimal representative of the inverse map$\,${\rm ;} moreover, any of them is an element of the symbolic power $I^{(n-1)}$
of degree $(m-1)(n-1)-1$.
\end{itemize}
\end{Proposition}
\demo
(a) By \cite[Theorem 3.2]{AHA}, it suffices to prove that the dimension of the $k$-subalgebra of $R$ generated
by the minors has dimension $n$, i.e., that $I$ has maximal analytic spread.
The case where $m=n$ follows from Proposition~\ref{general_cod2} (iii).
Now assume that $m>n$.
Since $R/I$ is Cohen--Macaulay and satisfies $\mu_P(I)\leq \hht P$, for $\hht P\leq n-1$ (Proposition~\ref{alma_mater} (i)),
 the result follows from \cite[Theorem 4.3]{UV}.

(b) For $m=n$ there is nothing to prove regarding the fiber type property, while
the symmetric algebra is even a complete intersection.
Thus, assume that $m>n$.
In this case the result follows from \cite[Theorem 1.3]{MoreyUlrich}.
In addition, the defining ideal of the Rees algebra $\mathcal{R}(I)$ is
$(I_1(\XX\cdot B), I_n(B))$,
where $B$ denotes the Jacobian dual matrix of $\mathcal{L}$.

\smallskip

(c) Since $I$ is of fiber type, a weak Jacobian dual matrix of $I$ as in \cite{AHA} coincides with
the transpose $B^t$ of the matrix introduced in the previous item;
$B^t$ is an $(m-1)\times n$ matrix of linear forms in the $\YY$-variables, whose rank  over the special
fiber of $I$ is $n-1$.
By part (a) and \cite{AHA}, any $(n-1)\times n$ submatrix has rank $n-1$ and its $n$
(ordered, signed) maximal minors are the coordinates
of a representative of the inverse map; thus, there are ${{m-1}\choose {n-1}}$ such representatives.

By construction, the degree of any one of these representatives (i.e., of its coordinates as elements of the
special fiber) is exactly $n-1$.
It follows from  Proposition~\ref{inversionfactor_is_symbolic} that each such representative gives rise to a source inversion
factor that is an element of the symbolic power $I^{(n-1)}$ and has degree $(m-1)(n-1)-1$.
\qed

\medskip

One gets immediately the following preamble to the subsequent main results.

\begin{Proposition}\label{general_symbolic_cod2}
Let $\mathcal{L}$ denote an $m\times (m-1)$
general linear matrix over $R=k[X_1,\ldots, X_n]$, with $m\geq n\geq 3$.
Set $I=I_{m-1}(\mathcal{L})\subset R$.
Then $I^{(r)}=I^r$ for $1\leq r\leq n-2$, and $D_j\in I^{(n-1)}\setminus I^{n-1}$,
where $D_j$ {\rm (}$j=1,\ldots, {{m-1}\choose {n-1}}${\rm )} are the source inversion
factors associated to a complete set of minimal representatives of the inverse map.
\end{Proposition}
\demo
The first assertion follows immediately from Theorem~\ref{associated_primes} and the second
assertion stems from Proposition~\ref{inversionfactor_is_symbolic}.
\qed

\smallskip

Thus far, the available features of the theory work for $m\geq n$.
In the subsequent part we come to grips with a richer amount of information, by focusing on  the cases where $m=n$ or $m=n+1$.
We will have to go a long way to obtain the nature of the corresponding symbolic Rees algebras.
Structure theorems for $m\geq n+2$ are this far unknown (see Remark~\ref{erratic}).

\subsection{The  symbolic algebra: Cremona case $m=n$}

The classical theory of plane Cremona maps in characteristic zero relates the Jacobian of a homaloidal net with the
principal curves of the corresponding Cremona map.
Our first proposition for this part is a far-fetched analogue of this result.

\begin{Proposition}\label{inversion_factor_is_determinant}{\rm (char$(k)=0$)}
Let $R=k[X_{1},\ldots, X_{n}]$ be a polynomial ring over a field $k$ of characteristic zero, with its standard grading
and let $\mathcal{L}=(\ell_{ij})$ be an $n\times (n-1)$ matrix whose entries are linear forms in $R$.
For every $i=1,\ldots, n$ write $\Delta_{i}$ for the signed $(n-1)$-minor of $\mathcal{L}$ obtained by omitting
the $i$-th row and let $\Theta=\Theta(\mathbf{\Delta})$ denote the Jacobian matrix of $\mathbf{\Delta}:=
\{\Delta_{1},\ldots,\Delta_{n}\}$.

If the ideal $I_{n-1}(\mathcal{L}):=(\mathbf{\Delta})\subset R$ is of linear type
then the rational map $\pp^{n-1}\dasharrow \pp^{n-1}$ defined by $\mathbf{\Delta}$
is a Cremona map and the associated source inversion factor is $\frac{1}{n-1}\det (\Theta)$.
\end{Proposition}
\demo
The first assertion to the effect that the map is birational is \cite[Examples 2.4]{cremona}
(also \cite[Theorem 3.12]{bir2003}).

We proceed to determine the source inversion factor.
Consider the Jacobian dual matrix of \cite{bir2003} which is the Jacobian matrix
with respect to $X_{1},\ldots, X_{n}$ of the linear forms in the target variables
$Y_1,\ldots, Y_n$ induced by the columns of $\mathcal{L}$.
This is the following matrix:

%\vspace{-10pt}

$$\left(\begin{array}{ccc}
\sum_{r=1}^{n}\frac{\partial \ell_{r1}}{\partial X_{1}}Y_{r}&\dots & \sum_{r=1}^{n}\frac{\partial \ell_{r1}}{\partial X_{n}}Y_{r}\\
\vdots&&\vdots\\
\sum_{r=1}^{n}\frac{\partial \ell_{rn-1}}{\partial X_{1}}Y_{r}&\ldots &\sum_{r=1}^{n}\frac{\partial \ell_{rn-1}}{\partial X_{n}}Y_{r}
\end{array}\right)$$

Now, by \cite[Examples 2.4]{cremona} the inverse map is defined by the (signed) maximal minors of this matrix.
Therefore, letting $\mathfrak{d}_{i}$ denote the signed $(n-1)$-minor of this matrix omitting the $i$th row,
by definition of the source inversion factor we are to show that the outcome of
evaluating $\mathfrak{d}_{i}$ via the map $Y_i\mapsto \Delta_i$ is $\frac{1}{n-1}\det (\Theta)\, X_i$.

To this purpose, we first note the following equality, where now $\Delta_i$ denotes the respective non-signed minor:
\vspace{-10pt}
{\small
$$\sum_{r=1}^{n}(-1)^{n+r}\frac{\partial \ell_{rj}}{\partial X_{k}}\Delta_{r}=\sum (-1)^{n+r+1} \ell_{rj}
\frac{\partial \Delta_{r}}{\partial X_{k}}, \; {\rm for}\; 1\leq k\leq n, \,1\leq j\leq n-1,$$
}
from which we gather:
{\small
$$\mathfrak{d}_{i}(\bf\Delta)=\det\left(\begin{array}{cccccc}
\sum_{r=1}^{n}(-1)^{n+r+1}\ell_{r1}\left(\begin{array}{c}
\frac{\partial\Delta_{r}}{\partial X_{1}}\\
\vdots\\
\frac{\partial\Delta_{r}}{\partial X_{i-1}}\\
\frac{\partial\Delta_{r}}{\partial X_{i+1}}\\
\vdots\\
\frac{\partial\Delta_{r}}{\partial X_{n}}\\
\end{array}\right)&\ldots&\sum_{r=1}^{n}(-1)^{m+r+1}\ell_{rn-1}\left(\begin{array}{c}
\frac{\partial\Delta_{r}}{\partial X_{1}}\\
\vdots\\
\frac{\partial\Delta_{r}}{\partial X_{r-1}}\\
\frac{\partial\Delta_{i}}{\partial X_{r+1}}\\
\vdots\\
\frac{\partial\Delta_{r}}{\partial X_{n}}\\
\end{array}\right)
\end{array}\right)$$
}

Write $[r_{1}\ldots r_{n-1}]$ for the $(n-1)$-minor of $\mathcal{L}$ with rows
$r_{1},\ldots,r_{n-1}$ and let $\alpha_{r_{1}\ldots r_{n-1}}:=(n+1)(n-1)+\sum_{s=1}^{n} r_{s}$.
By the multi-linearity of determinants, the result of evaluating $\mathfrak{d}_i$ is then

{\small
\begin{eqnarray}
\sum_{1\leq r_{1}<\ldots<r_{n-1}\leq n}&&\kern-10pt\left((-1)^{\alpha_{r_{1}\ldots r_{n-1}}}
\sum_{\sigma}(-1)^{\sigma} \ell_{\sigma(r_{1})1}\cdots \ell_{\sigma(r_{n-1})n-1}\right)\det\left(\begin{array}{ccccc}
\frac{\partial\Delta_{r_{1}}}{\partial X_{1}}&\ldots&\frac{\partial\Delta_{r_{n-1}}}{\partial X_{1}}\\
\vdots&\ddots&\vdots\\
\frac{\partial\Delta_{r_{1}}}{\partial X_{i-1}}&\ldots&\frac{\partial\Delta_{r_{n-1}}}{\partial X_{i-1}}\\
\frac{\partial\Delta_{r_{1}}}{\partial X_{i+1}}&\ldots&\frac{\partial\Delta_{r_{n-1}}}{\partial X_{i+1}}\\
\vdots&\ddots&\vdots\\
\frac{\partial\Delta_{r_{1}}}{\partial X_{n}}&\ldots&\frac{\partial\Delta_{r_{n-1}}}{\partial X_{n}}
\end{array}\right)\nonumber\\[5pt]
&=&\sum_{1\leq r_{1}<\ldots<r_{n-1}\leq n}(-1)^{\alpha_{r_{1}\ldots r_{n-1}}}[r_{1}\ldots r_{n-1}]\det\left(\begin{array}{ccccc}
\frac{\partial\Delta_{r_{1}}}{\partial X_{1}}&\ldots&\frac{\partial\Delta_{r_{n-1}}}{\partial X_{1}}\\
\vdots&\ddots&\vdots\\
\frac{\partial\Delta_{r_{1}}}{\partial X_{i-1}}&\ldots&\frac{\partial\Delta_{r_{n-1}}}{\partial X_{i-1}}\\
\frac{\partial\Delta_{r_{1}}}{\partial X_{i+1}}&\ldots&\frac{\partial\Delta_{r_{n-1}}}{\partial X_{i+1}}\\
\vdots&\ddots&\vdots\\
\frac{\partial\Delta_{r_{1}}}{\partial X_{n}}&\ldots&\frac{\partial\Delta_{r_{n-1}}}{\partial X_{n}}
\end{array}\right)\nonumber\\
&=&\det\left(\begin{array}{ccccc}
\frac{\partial\Delta_{1}}{\partial X_{1}}&\ldots&\frac{\partial\Delta_{n}}{\partial X_{1}}\\
\vdots&\ddots&\vdots\\
\frac{\partial\Delta_{1}}{\partial X_{i-1}}&\ldots&\frac{\partial\Delta_{n}}{\partial X_{i-1}}\\
\Delta_{1}&\ldots&\Delta_{n}\\
\frac{\partial\Delta_{1}}{\partial X_{i+1}}&\ldots&\frac{\partial\Delta_{n}}{\partial X_{i+1}}\\
\vdots&\ddots&\vdots\\
\frac{\partial\Delta_{1}}{\partial X_{n}}&\ldots&\frac{\partial\Delta_{n}}{\partial X_{n}}
\end{array}\right)=\frac{X_{i}}{n-1}\det\Theta\nonumber
\end{eqnarray}}

\noindent where we have expanded the determinant by Laplace  according to the $i$th row and used Euler's formula.

\begin{Corollary}\label{general_equal_case}{\rm (char$(k)=0$)}
Let $R=k[X_{1},\ldots, X_{n}]$ be a polynomial ring over a field $k$ of characteristic zero, with its standard grading
and let $\mathcal{L}$ be an $n\times (n-1)$ general linear matrix.
Then $I_{n-1}(\mathcal{L})$ is the base ideal of a Cremona map of $\pp^{n-1}$ and the associated
source inversion factor is $\frac{1}{n-1}\det (\Theta)$, where $\Theta$ denotes
the Jacobian matrix of the $(n-1)$-minors of $\mathcal{L}$.
\end{Corollary}
\demo The first assertion follows from Proposition~\ref{general_cod2}.
The second assertion is a consequence of
Proposition~\ref{inversion_factor_is_determinant}.
\qed

\begin{Remark}\rm An alternative to prove the second assertion of the previous corollary would come out
of Proposition~\ref{jac_vs_factor} by noticing that the inversion factor and $\det(\Theta)$ have the same degree.
Then it would suffice to argue that the latter is an irreducible
polynomial since the $(n-1)$-minors of $\mathcal{L}$ are sufficiently
general forms.
\end{Remark}

Here is the main theorem in the case $m=n$:

\begin{Theorem}\label{symbolic_m=n}%{\rm (char$(k)=0$)}
Let $\mathcal{L}$ denote an $n\times (n-1)$
general linear matrix over $R=k[X_1,\ldots, X_n]$, with $ n\geq 3$.
Set $I:=I_{n-1}(\mathcal{L})\subset R$
and let $\mathcal{R}^{(I)}$  denote its symbolic Rees algebra.
Then
\begin{enumerate}
\item[{\rm (a)}]  $\mathcal{R}^{(I)}$ is a Gorenstein normal domain.
\item[{\rm (b)}]  {\rm (char}$(k)=0${\rm )} $\mathcal{R}^{(I)}$ is generated by the $(n-1)$-minors of $\mathcal{L}$,
viewed in degree $1$ and by the source inversion factor of the Cremona map defined by
these minors, viewed in degree $n-1$.

Moreover, this inversion factor coincides with a nonzero scalar multiple of the Jacobian determinant of
the very minors.
\end{enumerate}

\end{Theorem}
\demo
(a) The symbolic Rees algebra
$\mathcal{R}^{(I)}$ of $I$ is a Gorenstein ring; indeed, it is a quasi-Gorenstein Krull
domain since $\hht I=2$ (\cite{ST}).
On the other hand, by the proof of \cite[Corollary 2.4 (b)]{dual}, one has an isomorphism
$\mathcal{R}^{(I)}\simeq \mathcal{R}(I)[t^{-1}]=R[It,t^{-1}]$, hence $\mathcal{R}^{(I)}$ is
finitely generated. Moreover, the latter is Cohen--Macaulay since
$\mathcal{R}(I)$ is Cohen--Macaulay by Proposition~\ref{m_greater_than_n_is _birational} (b).
It follows that $\mathcal{R}^{(I)}$ is a Gorenstein normal domain.

(b) To get the explicit generation, let $\mathfrak{d}_1,\ldots,\mathfrak{d}_n\in k[\YY]$ be forms
of the same degree, with $\gcd =1$, defining the inverse map and let $D\in R$ denote the corresponding
source inversion factor.
Write $J=(\mathfrak{d}_1,\ldots,\mathfrak{d}_n)\subset k[\YY]$.
By definition, one has
$$D=\mathfrak{d}_i(\Delta_{1},\ldots,\Delta_{n})/X_i,\, 1\leq i\leq n,$$
where $\mathbf{\Delta}:=
\{\Delta_{1},\ldots,\Delta_{n}\}$ are the (signed) minors generating $I$.
Identifying the two Rees algebras $\mathcal{R}_R(I)=R[It]\subset R[t]$ and
$\mathcal{R}_{k[\YY]}(J)=k[\YY][Ju]\subset k[\YY][u]$ by a $k$-isomorphism
that maps $Y_i\mapsto \Delta_i t$ and $X_i\mapsto \mathfrak{d}_iu$, then $D$ is identified with
$\mathfrak{d}_1/X_1$ in the common field of fractions.
Drawing on Proposition~\ref{alma_mater} (ii) (here we need char$(k)=0$), then the symbolic algebra is generated
by $It$ and $Dt^{n-1}$ as a consequence of
\cite[Corollary 7.4.3 (b)]{Wolmbook} (note that the notation for the two ideals is reversed
in the latter).

The additional statement follows from Corollary~\ref{general_equal_case} (again in characteristic zero).
\qed

\medskip

As an application of the results so far in the case $m=n$, we give an affirmative solution, in characteristic zero,
of the following conjecture stated in \cite[Section 2]{EiMa}:

\smallskip

\noindent{\bf Conjecture.}
If $I$ is the ideal of minors of a "generic" {\rm (}that is, random{\rm )}
$2 \times 3$ matrix of linear forms in $3$ variables, then the annihilator of $I^{(d)}/I^d$ is $F_1(I)^e$,
where $e$ is the greatest integer $\leq d/2$.

\smallskip

There is certainly a misprint in this statement since by definition the Fitting ideal $F_1(I)^e$ is the ideal of $2$-minors of
the matrix, which is the ideal $I$ itself.
The correct Fitting should be $F_2(I)^e$, the ideal of $1$-minors of the matrix.
But since the entries are general linear forms, this ideal is the irrelevant ideal $\fm:=(x,y,z)\subset R:=k[x,y,z]$
with $k$ a field.

\smallskip

\noindent{\bf Proof of the conjecture.}
%More is actually true but this will suffice.
%Since $\deg(D)=3$, the inclusions $xD\in I^2$, $yD\in I^2$ and $zD\in I^2$ are $k$-linear combinations.

Let $\phi$ denote the given matrix.
A consequence of Theorem~\ref{symbolic_m=n} (b) above is that $I^{(2)}=(I^2, D)$, where $D\in R$ is the
inversion factor of the Cremona map defined by the $2$-minors of $\phi$, and, moreover,
for every $d\geq 1$ the following equalities hold
$$I^{(d)}=
\left\{
\begin{array}{ll}
(I^{(2)})^{\frac{d}{2}} &\mbox{if $d$ is even}\\
I\,(I^{(2)})^{\frac{d-1}{2}} &\mbox{if $d$ is odd}
\end{array}
\right.
$$

By definition of the inversion factor, one has $D\fm\in I^2$.
It follows that the annihilator of $I^{(2)}/I^2$ is $\fm$.
We also know that $\deg(D)=3$ since the inverse map to the Cremona map defined by the $2$-minors
is also defined by forms of degree $2$, and so $\deg(D)=2.2-1=3$.

Consider separately the even and the odd cases.

\smallskip

$d$ {\sc even.}

\smallskip

One has $\fm^{d/2} I^{(d)}= \fm^{d/2}(I^{(2)})^{d/2}= (\fm I^{(2)})^{d/2}\subset (I^2)^{d/2}=I^d.$

Conversely, let $f\in \fm$ be a form such that $fI^{(d)}\subset I^d$.
Since the annihilator of $I^{(2)}/I^2$ is the entire maximal ideal $\fm$, it suffices to show that
$\deg(f)\geq d/2$.
Since $I^{(d)}= (I^{(2)})^{d/2}$ and $I^{(2)}=(I^2,D)$, in particular we get $fD^{d/2}\in I^d$.
Reading degrees on both sides one has that $\deg(f)+3d/2\geq 2d$. Therefore,
$\deg(f)\geq d/2$, hence $f\in \fm^{d/2}$ as required.

\medskip

$d$ {\sc odd.}

\smallskip

One has
$$\fm^{(d-1)/2} I^{(d)}= \fm^{(d-1)/2}\,I(I^{(2)})^{(d-1)/2}=I\, (\fm I^{(2)})^{(d-1)/2}\subset
I\,(I^2)^{(d-1)/2}=I\,I^{d-1}=I^d.$$

The hypothesis is that $fD^{(d-1)/2}\, I\subset I^d$.
In particular, taking a minor $\Delta$ among the generators of $I$, we find $f\,D^{(d-1)/}\,\Delta\subset I^d$.
Again, reading degrees, we get the inequality $\deg(f) +  3((d-1)/2) +2 \geq 2d$, from which follows that
$\deg(f)\geq (d-1)/2$.

We conclude as before.
Since in the odd case, $(d-1)/2=\floor{d/2}$, we are done.
\qed

\subsection{The  symbolic algebra: implicitization case $m=n+1$}

We will now assume that $m=n+1$.

\subsubsection{Homological prelims}

The arguments in this part will draw on the following results of independent interest.
To describe their contents, recall that $\mathcal{S}_{n-1}(I)\simeq I^{n-1}$ by Corollary~\ref{F_1}.
Therefore, by \cite{AVRAMOV, TCHERNEV, Wey} one has a free resolution of $I^{n-1}$
$$\mathcal{K}_{n-1}\;:\;0\rar F_{n-1}\rar F_{n-2}\rar\ldots\rar F_1\rar F_0\rar0$$
where $$F_i:=\bigwedge^{i}R^{n}\otimes_R \mathcal{S}_{(n-1)-i}(R^{n+1})$$
and $d:F_i\rar F_{i-1}$
is given by $$d(e_1\wedge\ldots\wedge e_i\otimes g):=\sum_{l=1}^{i}e_1\wedge\ldots\wedge\widehat{e_l}
\wedge\ldots\wedge e_i\otimes\phi(e_l)g,$$
with $\{e_1, \ldots, e_n\}$ denoting a basis of $R^n$ and  $\phi:R^n \rar R^{n+1}$
 standing for the map defined by the $(n+1)\times n$ presentation matrix  $\mathcal{L}=
(\ell_{ij})$ of the ideal $I$.

Consider the $R$-dual map to $d_{n-1}: F_{n-1}\rar F_{n-2}$.
Since $I^{n-1}$ is generated in (standard) degree $n(n-1)$, after identification and taking in account the degrees shift,
the dual map is of the form
 \begin{equation}\label{dual_degrees}
\eta:=d_{n-1}^*: R^N((n+1)(n-1)-1)\;{\rar}\; R^n((n+1)(n-1)),
 \end{equation}
where $N=(n+1){{n}\choose {2}}$.
Let $M$ denote the cokernel of $\eta$.
Shifting by $-((n+1)(n-1))$, we get a homogeneous presentation
 \begin{equation*}%\label{dual_degrees}
R^N(-1)\;\stackrel{\eta}{\rar}\; R^n\;\rar\; M(-(n+1)(n-1))\rar 0.
 \end{equation*}

\begin{Theorem}\label{funny_module}
With the above notation, there is a homogeneous isomorphism
$$M(-(n+1)(n-1))\simeq R^n/(\XX)R^n=k^n.$$
\end{Theorem}
\demo
Picking up from the above preliminaries, let us make explicit the dual map to $d_{n-1}: F_{n-1}\rar F_{n-2}$.
Note that
$$F_{n-1}=\bigwedge^{n-1}R^n\otimes_R S_0(R^{n+1})\simeq R^n, \;\; F_{n-2}=\bigwedge^{n-2}R^n\otimes_R S_1(R^{n+1})
\simeq R^{{n}\choose {n-2}}\otimes_R R^{n+1}.$$
Applying these identifications, the basis vector $e_1\wedge\cdots\wedge \widehat{e_k}\wedge\cdots\wedge e_n$ gets identified with $e_k$
and we write $a_{1,\ldots ,\hat{k},\ldots,\hat{l},\ldots, n}$ for a basis vector of $R^{{n}\choose {n-2}}$ corresponding to
$e_1\wedge\cdots\wedge \widehat{e_j}\wedge\cdots\wedge\widehat{e_l}\wedge\cdots\wedge e_n$.
Further, let $\{b_1,\ldots,b_{n+1}\}$ stand for a basis of $R^{n+1}$
With this notation, for $k=1,\ldots,n$, the map is quite simply

\begin{eqnarray*}
e_k&\mapsto &\sum_{l=1}^{n-1} a_{1,\ldots ,\hat{k},\ldots,\hat{l},\ldots, n}
\otimes\phi(e_l)=\sum_{l=1}^{n-1} a_{1,\ldots ,\hat{k},\ldots,\hat{l},\ldots, n}
\otimes\sum_{i=1}^{n+1} \ell_{il}b_i\\
&=& \sum_{i=1}^{n+1}\sum_{l=1}^{n-1} \ell_{il}\; a_{1,\ldots ,\hat{k},\ldots,\hat{l},\ldots, n}
\otimes b_i,
\end{eqnarray*}
where $\mathcal{L}=(\ell_{ij})$ is as above.

From this the transposed matrix has the following block shape
$$\eta=\left(M_{n-1,n}|\ldots |M_{1,n}|\ldots |M_{j-1,j}|\ldots|M_{1,j}|\ldots| M_{1,2}\right),$$
where, for $1\leq i\leq j\leq n,$ $M_{ij}$ is the following $n\times (n+1)$ matrix up to signs

 $$M_{i,j}=\left(\begin{array}{cccccccc}
 0&0&\ldots&0\\
 \vdots&\vdots&\ldots&\vdots\\
 \ell_{1i}&\ell_{2i}&\ldots&\ell_{(n+1)i}\\
  \vdots&\vdots&\ldots&\vdots\\
   \ell_{1j}&\ell_{2j}&\ldots&\ell_{(n+1)j}\\
    \vdots&\vdots&\ldots&\vdots\\
     0&0&\ldots&0\\
 \end{array}\right)
 \begin{array}{cc}
 &\\[3pt]
 \leftarrow\;(n+1-j)\mbox{th row}
 &\\
 &\\
 &\\
 \leftarrow\;(n+1-i)\mbox{th row}
 &\\
 &
 \end{array}$$
Next let $\widetilde{M_{i,j}}$ denote the submatrix of $M_{i,j}$ consisting of the first $n$ columns and consider
the following block submatrix of $\eta$
$$(\widetilde{M_{n-1,n}}|\ldots| \widetilde{M_{1,n}}\,|\,\widetilde{M_{n-2,n-1}}),$$
consisting of $n$ square blocks of order $n$ each; in particular, the matrix has $n^2$ columns.

We claim that the $R$-submodule of $R^n$ generated by the columns of the above matrix coincides with $(\XX)R^n$.
For this, since the columns have standard degree $1$, it suffices to show that the columns are $k$-linearly independent
as elements of the $k$-vector space $((\XX)R)_1$.

Suppose that a nontrivial $k$-linear combination of these columns vanishes, with coefficients $\alpha_1,\ldots,\alpha_{n^2}\in k$.
Grouping the coefficients corresponding to the variables $X_1,\ldots,X_n$ one gets an $n\times n$ linear system
with coefficients in $k$ such that $\{X_1,\ldots,X_n\}$ is a non zero solution.
But then every row of the system gives a $k$-linear relation of these variables.
Clearly this is only possible if all the coefficients of this system vanish.
Writing this condition as a new square linear system, this time around of order $n^2$ with solution $\{\alpha_1,\ldots,\alpha_{n^2}\}$
and appropriate coefficients in $k$.
Since the latter coefficients are nothing but the coefficients of all linear forms $\ell_{ij}$, they can be expressed as partial
derivatives of these forms, so the corresponding $n^2\times n^2$ matrix has the following form (up to signs)
$$\Theta=\left(\begin{array}{ccccccccccccccccccccc}
\Theta_{n-1}&\Theta_{n-2}&\Theta_{n-3}&\ldots&\Theta_2& \Theta_{1}&{\bf 0}\\
\Theta_n&{\bf 0}&{\bf 0}&\ldots&{\bf 0}&{\bf 0}&\Theta_{n-2}\\
{\bf 0}&\Theta_{n}&{\bf 0}&\ldots&{\bf 0}&{\bf 0}&\Theta_{n-1}\\
{\bf 0}&{\bf 0}&\Theta_{n}&\ldots&{\bf 0}&{\bf 0}&{\bf 0}\\
\vdots&\vdots&\vdots&\ldots&\vdots&\vdots&\vdots\\
{\bf 0}&{\bf 0}&{\bf 0}&\ldots&\Theta_n&{\bf 0}&{\bf 0}\\
{\bf 0}&{\bf 0}&{\bf 0}&\ldots&{\bf 0}&\Theta_n&{\bf 0}
\end{array}\right)$$
where $\Theta_i$ is the transpose of the Jacobian matrix of $\{\ell_{i1},\ldots, \ell_{in}\}$ and ${\bf 0}$ denotes the null
matrix of order $n$.
We note that $\Theta_i$ is non-singular since $\{\ell_{i1},\ldots, \ell_{in}\}$ is a set of $k$-linearly independent
$1$-forms.
The system has only the trivial solution if and only if the determinant of this matrix does not vanish.
One can see that, after appropriate elementary row operations, the above determinant is non-vanishing if and only if
the following matrix has nonzero determinant
$$\left(\begin{array}{ccccccccccccccccccccc}
I&{\bf 0}&{\bf 0}&\ldots&{\bf 0}&{\bf 0}&\Theta_{n-2}\\
{\bf 0}&I&{\bf 0}&\ldots&{\bf 0}&{\bf 0}&\Theta_{n-1}\\
{\bf 0}&{\bf 0}&I&\ldots&{\bf 0}&{\bf 0}&{\bf 0}\\
\vdots&\vdots&\vdots&\ldots&\vdots&\vdots&\vdots\\
{\bf 0}&{\bf 0}&{\bf 0}&\ldots&I&{\bf 0}&{\bf 0}\\
{\bf 0}&{\bf 0}&{\bf 0}&\ldots&{\bf 0}&I&{\bf 0}\\
{\bf 0}&{\bf 0}&{\bf 0}&\ldots&{\bf 0}&0&\Omega
\end{array}\right)$$
where $I$ stands for the $n\times n$ identity matrix and $\Omega=\Theta_{n-1}\Theta_{n-2}-\Theta_{n-2}\Theta_{n-1}$.
Thus, $\det(\Theta)\neq 0$ if and only if $\det(\Omega)\neq 0.$
Now, the entries of the matrices $\Theta_{n-1},\Theta_{n-2}$ are among the coefficients
of the entries of the matrix $\mathcal{L}=(\ell_{ij})$.
Since these are random, $\det(\Omega)\neq 0.$

\smallskip

Now, to conclude, we have shown that the image of the map $\eta$ in (\ref{dual_degrees}) is the $R$-submodule $(\XX)R^n$.
Therefore, $M(-(n+1)(n-1))\simeq R^n/(\XX)R^n$ as required.
\qed
\begin{Example}\label{Tcher}\rm The above discussion has many common points with \cite[Section 8.2]{Wolmbook} which
treats the case of linearly presented perfect ideals in dimension $n=3$.
However, the above proof draws on the hypothesis that $\mathcal{L}$ is a
general linear matrix -- and, in fact, it may be false for other linearly presented ideals.
We are indebted to A. Tchernev for having provided us with the following counter-example to
Theorem~\ref{funny_module}
in the context of arbitrary linearly presented ideals:
\begin{eqnarray}\label{tcher2}
\phi=
\left(
\begin{array}{ccc}
X_1  &   X_2     &  X_3 \\
X_2  &  X_3    &   0  \\
X_3  &  0    &  X_1 \\
0    &  X_1    &  X_2
\end{array}
\right)
\end{eqnarray}
Here the vector space dimension of the linear forms in Im$(\eta)$ is $8$, where $\eta$ denotes the correspondingly
defined matrix as in the proposition.
We note that by changing $6$ out of the $9$ nonzero entries of the Tchernev matrix into general $1$-forms,
the resulting matrix gives the maximal value $9$ for the vector space dimension of the linear forms in Im$(\eta)$.
\end{Example}

\begin{Proposition}\label{local_duality}
Let $\mathcal{L}$ denote an $(n+1)\times n$
general linear matrix over $R=k[X_1,\ldots, X_n]$, with $n\geq 3$.
Set $I=I_{n}(\mathcal{L})\subset R$.
Then
 $$I^{(n-1)}/I^{n-1}\simeq k^n (-(n(n-1)-1)),$$
 as graded $R$-modules.
\end{Proposition}
\demo
By Theorem~\ref{funny_module}, one has a (shifted) homogeneous isomorphism
$$M(-n)\simeq k^n (n(n-1)-1).$$
On the other hand, by definition there is a homogeneous isomorphism
$$M\simeq {\rm Ext}_R^n(R/I^{n-1}, R).$$
Therefore, it obtains
\begin{eqnarray*}
I^{(n-1)}/I^{n-1}&\simeq & H_{(\XX)}^0(R/I^{n-1}) \quad (\mbox{since $I^{(n-1)}/I^{n-1}$ has finite length})\\
&\simeq & {\rm Hom}_R({\rm Ext}_R^n(R/I^{n-1}, R(-n)), E(k)) \quad (\mbox{by graded local duality})\\
&\simeq & {\rm Hom}_R(M(-n), E(k))\simeq {\rm Hom}_R(k^n(n(n-1)-1), E(k))\\
&\simeq & {\rm Hom}_R(k, E(k))^{\oplus n} (-(n(n-1)-1))\simeq k^n (-(n(n-1)-1)),
\end{eqnarray*}
where the last isomorphism is given in \cite[Lemma 3.2.7 (b)]{BHbook}.
\qed

\begin{Example}\label{symbolic_cyclic} \rm Corollary~\ref{local_duality} fails for arbitrary perfect ideals of codimension $2$
admitting linear presentation. For $n=3$, Example~\ref{Tcher} is a counter-example. Letting $I\subset R=k[X_1,X_2,X_3]$
denote the ideal of $3$-minors, then $I^{(2)}/I^2$ is a cyclic $R$-module generated by the residue class of a form $F\in I^{(2)}$
of degree $4<n(n-1)-1=5$.
The map defined by the minors is still birational onto the image, with inversion factors $X_1F,X_2F,X_3F$.
In particular, the latter are not minimal generators of $I^{(2)}$.
Even if we slightly `perturb' Tchernev's matrix the result equally fails, such as in the following matrix
$$
\left(\begin{array}{ccc}
X_1&X_2&X_3\\
X_2&X_3&0\\
X_3&0&X_1-X_2\\
0&X_1-X_3&X_2-X_3
\end{array}
\right).
$$
Perturbing more entries, such as in the following matrix
$$
\left(\begin{array}{ccc}
X_1&X_2&X_3\\
X_2&X_3&X_1-X_2\\
X_3&X_1-X_3&X_2-X_3\\
X_1+X_2&X_2+X_3&X_1+X_3
\end{array}
\right),
$$
the result of Proposition~\ref{local_duality} still holds true -- and so does the one of Theorem~\ref{funny_module}.
However, now the statement in Theorem~\ref{The D's} (i) fails as those generators have a nontrivial common divisor.
Cooking up some of these examples require some extra care  to make sure that $I$ is a radical ideal, otherwise the whole
known repository of symbolic power theory crumbles down.
Therefore, slightly perturbing a linear $(n+1)\times n$ matrix whose $n$-minors generate
a radical ideal may lead us astray.
As an example, changing the lower right corner entry of the first of the above matrices into $X_1-X_2+X_3$
gives a non-radical ideal.
\end{Example}

\subsubsection{The trick of the transposed Jacobian dual}

A good deal of the subsequent development rests on a simple construction.

Namely, let ${\bf \Delta}=\{\Delta,\ldots, \Delta_{n+1}\}$ denote the signed maximal minors of $\mathcal{L}$.
Let $B$ denote the Jacobian dual matrix of $\mathcal{L}$, whose entries belong to the polynomial ring
$k[\YY]=k[Y_1,\ldots,Y_n,Y_{n+1}]$. By definition, one has an equality $\YY\cdot\mathcal{L}=\XX\cdot B^t$,
where the superscript $^t$ denotes transpose.
We can similarly write an equality $\YY\cdot\mathcal{L}'=\ZZ\cdot B$, for a unique matrix
$\mathcal{L}'$ whose entries are linear forms in a set of duplicate variables $\ZZ$ of $\XX$.

We observe  that $\mathcal{L}'$ only differs from $\mathcal{L}$  by the rearrangement of the (same)
coefficients of the linear forms. Since the notion of general linear forms is dictated by the randomness of the
total set of coefficients, it follows that $\mathcal{L}'$ it too is a matrix whose entries are general linear
forms in the variables $\ZZ$.
Therefore, by Proposition~\ref{m_greater_than_n_is _birational} (a), its $n$-minors
${\boldsymbol\delta}=\{\delta_1,\ldots,\delta_{n+1}\}$ define a birational map onto the image,
with $B^t$ as its Jacobian dual matrix and corresponding set $\{d_1(\ZZ), \ldots, d_n(\ZZ)\}$ of source inversion
factors associated to a complete set of minimal representatives of the corresponding inverse map.

As usual, the set of $n$-minors is taken with the correct signs.
Keeping the above notation, one has the following basic structural result:

\begin{Theorem}\label{The D's}{\rm (char$(k)=0$)}
Let $\mathcal{L}$ denote an $(n+1)\times n\; (n\geq 3)$
general linear matrix over $R=k[\XX]=k[X_1,\ldots, X_n]$ with $n$-minors
${\boldsymbol\Delta}=\{\Delta_1,\ldots,\Delta_{n+1}\}$,
 let $\{D_1, \ldots, D_n\}$ be as in {\rm Proposition~\ref{general_symbolic_cod2}},
 and $\{d_1(\ZZ), \ldots, d_n(\ZZ)\}$ as above.
Then:
\begin{enumerate}
\item[{\rm (i)}] $\{D_1, \ldots, D_n\}\subset R$ and $\{d_1(\ZZ), \ldots, d_n(\ZZ)\}\subset k[\ZZ]$
both generate ideals of codimension $2$.
\item[{\rm (ii)}] $\{D_1, \ldots, D_n\}$ defines a Cremona map
$\mathfrak{D}$ of $\pp^{n-1}$ whose inverse map is $(d_1(\ZZ): \cdots : d_n(\ZZ))$.
\item[{\rm (iii)}] Writing $I:=(\Delta_1,\ldots,\Delta_{n+1})$, the $R$-module $I^{(n-1)}/I^{n-1}$ is minimally generated by the classes of
$D_1,\ldots, D_n${\rm ;} in particular, the symbolic power $I^{(n-1)}$ is minimally generated by $D_1,\ldots, D_n$
and by the minimal generators of $I^{n-1}$ which are not of the form $X_iD_j, 1\leq i,j\leq n$.
\item[{\rm (iv)}]  The source inversion factor of $\,\mathfrak{D}$ is the $(n-1)$th power of an
element $E\in I^{(n(n-1)-1)}$.

\noindent {\sc Supplement}: Moreover, $E$ coincides with the Jacobian determinant of $D_1,\ldots, D_n$
whenever the latter is irreducible.

\item[{\rm (v)}] The minimal graded resolution of the ideal $(D_1, \ldots, D_n)\subset R$ is
\begin{equation}\label{resolucaodosDs}
0\rar\; R(-n^2)\;\stackrel{\XX^t}{\rar}\; R(-(n^2-1))^n\;\stackrel{\Psi}{\rar}R(-(n(n-1)-1))^{n}\;\rar R,
\end{equation}
where $\Psi$ denotes the Jacobian dual matrix of the signed $n$-minors of $\mathcal{L}$ evaluated orderly on
these signed minors, while $\XX^t$ stands for the transpose of the vector of the source variables.
\end{enumerate}
\end{Theorem}
\demo
(i) We only discuss the ideal $(D_1,\ldots,D_n)$ since the line of argument is analogous for
$(d_1(\ZZ), \ldots, d_n(\ZZ))$.

Being a subideal of $I:=(\Delta_1,\ldots,\Delta_{n+1})$, the codimension of $(D_1,\ldots,D_n)$ is at most $2$.
Thus, it suffices to show that it is at least $2$.
Start from scratch by observing that $k[{\boldsymbol\Delta}]\simeq k[\YY]/(\boldsymbol\beta)$, where $\boldsymbol\beta:=\det (B)$
and $B$ stands for the Jacobian dual matrix of ${\boldsymbol\Delta}$.
Since $\boldsymbol\beta({\boldsymbol\Delta})=0$, the chain rule of derivatives gives the short polarization  complex
\begin{equation}\label{polarization_complex}
R\stackrel{\partial}{\lar} R^{n+1}\stackrel{\Theta}{\lar} R^n,
\end{equation}
where $\Theta$ denotes the transposed Jacobian matrix of ${\boldsymbol\Delta}$ and $\partial$ is the transpose of
$$\left[\frac{\partial \boldsymbol\beta}{\partial Y_1}({\boldsymbol\Delta})\,\ldots \,
\frac{\partial \boldsymbol\beta}{\partial Y_{n+1}}({\boldsymbol\Delta})\right].$$
On the other hand, since $\dim k[{\boldsymbol\Delta}]=n$, the rank of $\Theta$ is $n$ (since char$(k)=0$), hence $\ker(\Theta)$ is generated
by the single (column) vector whose $j$th coordinate
is the $n$-minor of $\Theta$ omitting the $j$th column of $\Theta$ further divided by the $\gcd$ of all the $n$-minors.
Since ${\boldsymbol\Delta}$ are maximal minors of a general linear matrix, they are sufficiently general $n$-forms,
and so are any of their derivatives (the entries of $\Theta$).
Since having a proper common divisor is a closed condition on the coefficients while the coefficients
of the entries are products and sums of random coefficients, then the ideal $I_n(\Theta)$ generated
by the maximal minors of $\Theta$ has codimension $2$.

This implies that $\ker(\Theta)$ is generated by a single vector in degree $(n-1)n$ (the degree of an $n$-minor of $\Theta$).
On the other hand, a simple calculation shows that the coordinates of $\partial$ are also of degree $n(n-1)$.
Since by (\ref{polarization_complex}) the $j$th coordinate of $\partial$ is a multiple of the $n$-minor of $\Theta$ omitting the $j$th column,
we must conclude that the ideals $\left(\frac{\partial \boldsymbol\beta}{\partial Y_1}({\boldsymbol\Delta}),\,\ldots ,\,
\frac{\partial \boldsymbol\beta}{\partial Y_{n+1}}({\boldsymbol\Delta})\right)$ and $I_n(\Theta)$ coincide.

In particular, the first of these ideals has codimension $2$.
We proceed to show that it is further contained in the ideal $(D_1,\ldots, D_n)$, thus showing that the latter
has  codimension at least $2$.

Let $\mathcal{L}=(\ell_{ij})$ denote the given general linear $(n+1)\times n$ matrix.
Then
$$B^t=\left(\begin{array}{cccccccc}
\sum_{r=1}^{n+1}\frac{\partial \ell_{r,1}}{\partial X_1}Y_r&\sum_{r=1}^{n+1}\frac{\partial \ell_{r,2}}{\partial X_1}Y_r&
\ldots&\sum_{r=1}^{n+1}\frac{\partial \ell_{r,n}}{\partial X_1}Y_r\\[7pt]
\sum_{r=1}^{n+1}\frac{\partial \ell_{r,1}}{\partial X_2}Y_r&\sum_{r=1}^{n+1}
\frac{\partial \ell_{r,2}}{\partial X_2}Y_r&\ldots&\sum_{r=1}^{n+1}\frac{\partial \ell_{r,n}}{\partial X_2}Y_r\\[5pt]
\vdots&\vdots&\ddots&\vdots\\[5pt]
\sum_{r=1}^{n+1}\frac{\partial \ell_{r,1}}{\partial X_n}Y_r&\sum_{r=1}^{n+1}
\frac{\partial \ell_{r,2}}{\partial X_n}Y_r&\ldots&\sum_{r=1}^{n+1}\frac{\partial \ell_{r,n}}{\partial X_n}Y_r
\end{array}\right).$$
Expanding the determinant of $B$,  it obtains (up to signs)
$$\boldsymbol\beta=\sum_{1\leq j_1\leq\cdots\leq j_n\leq n}\left[\left(\sum_{r=1}^{n+1}\frac{\partial \ell_{r,j_1}}{\partial X_1}Y_r\right)\cdots
\left(\sum_{r=1}^{n+1}\frac{\partial \ell_{r,j_n}}{\partial X_1}Y_r\right)\right].$$
Taking the $k$th derivative yields
\begin{eqnarray*}\frac{\partial \boldsymbol\beta}{\partial X_k}&=&\sum_{1\leq j_1\leq\cdots\leq j_n\leq n}
\left[\sum_{1\leq s\leq n}\left(\sum_{r=1}^{n+1}\frac{\partial \ell_{r,j_1}}{\partial X_1}Y_r\right)\cdots
\left(\frac{\partial \ell_{k,j_s}}{\partial X_s}\right)\cdots
\left(\sum_{r=1}^{n+1}\frac{\partial \ell_{r,j_n}}{\partial X_1}Y_r\right)\right]\\
&=& \sum_{1\leq s\leq n}
\left[ \sum_{1\leq j_1\leq\cdots\leq j_n\leq n}  \left(\sum_{r=1}^{n+1}\frac{\partial \ell_{r,j_1}}{\partial X_1}Y_r\right)\cdots
\left(\frac{\partial \ell_{k,j_s}}{\partial X_s}\right)\cdots
\left(\sum_{r=1}^{n+1}\frac{\partial \ell_{r,j_n}}{\partial X_1}Y_r\right)\right]
\end{eqnarray*}
Note that for any given $1\leq s\leq n$, the expression inside the square brackets in the last line above is
(up to signs) the determinant of the matrix

$$B_s:=\left(\begin{array}{cccccccc}
\sum_{r=1}^{n+1}\frac{\partial \ell_{r,1}}{\partial X_1}Y_r&\sum_{r=1}^{n+1}\frac{\partial \ell_{r,2}}{\partial X_1}Y_r&
\ldots&\sum_{r=1}^{n+1}\frac{\partial \ell_{r,n}}{\partial X_1}Y_r\\[7pt]
\vdots&\vdots&\ddots&\vdots\\[5pt]
\frac{\partial \ell_{k,1}}{\partial X_s} & \frac{\partial \ell_{k,2}}{\partial X_s} & \ldots &
\frac{\partial \ell_{k,n}}{\partial X_s}\\[5pt]
\vdots&\vdots&\ddots&\vdots\\[5pt]
\sum_{r=1}^{n+1}\frac{\partial \ell_{r,1}}{\partial X_n}Y_r&\sum_{r=1}^{n+1}
\frac{\partial \ell_{r,2}}{\partial X_n}Y_r&\ldots&\sum_{r=1}^{n+1}\frac{\partial \ell_{r,n}}{\partial X_n}Y_r
\end{array}\right).$$

Expanding this determinant once again, this item around by Laplace along the $i$th row of $B_s$, gives
$\det(B_s)=\sum_{t=1}^n \frac{\partial \ell_{k,t}}{\partial X_s}\,\sigma_t^{[s]}$,
where $\sigma_t^{[s]}$ denotes the $(n-1)$-minor of $B_s$ omitting the $s$th row and the $t$th column.
Coming from the other end, for given $s$, $(\sigma_1^{[s]}:\cdots :\sigma_n^{[s]})$ is a representative of
the inverse map to the map defined by $\boldsymbol\Delta$ (\cite[Theorem 2.18, Supplement]{AHA}).
By definition, say, $D_s$ is the source inversion factor corresponding to this representative,
hence $\sigma_t^{[s]}(\boldsymbol\Delta)=X_tD_s$, for $1\leq s,t \leq n$.
Assembling the information, we get
\begin{eqnarray*}\frac{\partial \boldsymbol\beta}{\partial X_k}(\boldsymbol\Delta)&=&\det(B_1)(\boldsymbol\Delta)
 +\cdots + \det(B_n)(\boldsymbol\Delta)=
\sum_{t=1}^n \frac{\partial \ell_{k,t}}{\partial X_1}\sigma_t^{[1]}(\boldsymbol\Delta) +\cdots +
\sum_{t=1}^n \frac{\partial \ell_{k,t}}{\partial X_n}\sigma_t^{[n]}(\boldsymbol\Delta)\\
&=& \left(\sum_{t=1}^n \frac{\partial \ell_{k,t}}{\partial X_1}X_t\right) D_1 +\cdots +
\left(\sum_{t=1}^n \frac{\partial \ell_{k,t}}{\partial X_n}X_t\right) D_n,
\end{eqnarray*}
which proves our contention.

\smallskip

(ii)
Let ${\boldsymbol\delta}=\{\delta_1,\ldots,\delta_{n+1}\}\subset k[\ZZ]$ stand for the $n$-minors of the
general linear matrix  $\mathcal{L}'$ as explained above. We have seen in the preliminaries of this section that they
define a birational map onto the image, with $B^t$ as its Jacobian dual matrix.
Thus, for any $j\in \{1,\ldots, n\}$, the $(i,j)$-cofactors of $B^t$
$$\{B^{t}_{j1},\ldots,B^{t}_{jn}\}$$
taken modulo $\det(B)$ define an inverse to the map defined by ${\boldsymbol\delta}$.
By Lemma\ref{uniqueness_of_target_factor} this yields the following structural
congruencies
 \begin{equation}\label{igualdadeestrutural}
\delta_i(B^t_{j1},\ldots, B^t_{jn})\equiv E_{j}Y_i\mod( \det B^t),
\end{equation}
where $E_j$ denote the corresponding target inversion factor.

{\bf Claim.} $I_1(\XX\cdot B^t(\boldsymbol\delta))$ is contained in the presentation ideal of the
Rees algebra of the ideal $(D_1, \ldots ,D_n)\subset k[\XX]$, defined over the ring $k[\XX,\ZZ]$.

To see this it suffices to prove that the entries of $\XX\cdot B^t$ vanish by evaluating $Y_k\mapsto \delta_k(D_1,\ldots,D_n)$,
$k=1,\ldots,n+1$, or, equivalently, by evaluating $Y_k\mapsto \delta_k(X_nD_1,\ldots,X_nD_n)$,
$k=1,\ldots,n+1$.
Letting, as previously, $\boldsymbol\Delta=\{\Delta_1,\ldots, \Delta_{n+1}\}$ denote the signed $n$-minors of
$\mathcal{L}$, one has the relations
\begin{equation}\label{transposta}
X_nD_i=B_{in}(\Delta_1,\ldots,\Delta_{n+1})=B^t_{ni}(\Delta_1,\ldots,\Delta_{n+1}),
\end{equation}
since $D_i$ is inversion factor for $\boldsymbol\Delta$, where $B_{ij}$ is the cofactor of $B$ corresponding to the entry indexed by $(i,j)$
and $B_{in}(\Delta_1,\ldots,\Delta_{n+1})$ is the result of evaluating this cofactor on $\boldsymbol\Delta$.
Since $B_{ij}=B_{ji}^t$, one gets
{\small
\begin{eqnarray}
\XX\cdot B^t(\delta(X_nD_1,\ldots,X_nD_n))&=&\XX\cdot\left(\begin{array}{ccccccc}
\sum_{i=1}^{n+1}\frac{\partial \ell_{i1}}{\partial X_1 }\delta_i(X_n{\bf D})&
\ldots&\sum_{i=1}^{n+1}\frac{\partial\ell_{in}}{\partial X_1 }\delta_i(X_n{\bf D})\\
\vdots&\ddots&\vdots\\
\sum_{i=1}^{n+1}\frac{\partial\ell_{i1}}{\partial X_3}\delta_i(X_n{\bf D})&
\ldots&\sum_{i=1}^{n+1}\frac{\partial\ell_{in}}{\partial X_n}\delta_i(X_n{\bf D})
\end{array}\right)\\
&\vspace{0.5cm}&\nonumber\\\label{eqrees1}
&=&\left(\sum_{i=1}^{n+1}\ell_{i1}\delta_i(X_n{\bf D}),\ldots,\sum_{i=1}^{n+1}\ell_{in}\delta_i(X_n{\bf D})\right)\\
&\vspace{0.5cm}&\nonumber\\\label{eqrees2}
&=&E_n({\bf \Delta})\left(\sum\ell_{i1}\Delta_i,\ldots,\sum\ell_{in}\Delta_i\right)\\\label{eqrees3}
&=&(0,\ldots,0)
\end{eqnarray}}
where equality (\ref{eqrees2}) follows from (\ref{igualdadeestrutural}), (\ref{transposta}) and (\ref{eqrees1}) -- keep
 in mind that the result of evaluating $\det(B^t)$ by $Y_K\mapsto \Delta_k$ is zero. As to equality (\ref{eqrees3}), it is a
 consequence of (\ref{eqrees2}) using that $\mathcal{L}$ is a syzygy matrix of ${\bf\Delta}.$
This proves the claim.

As a consequence, the matrix $B({\boldsymbol\delta})$ is a submatrix of the full Jacobian dual matrix of
${\bf D}:=\{D_1,\ldots, D_n\}$.
On the other hand, we have $\det (B({\boldsymbol\delta}))= (\det (B))({\boldsymbol\delta})=(\det (B^t))({\boldsymbol\delta})=0$
since $\det (B^t)$ is a polynomial relation of ${\boldsymbol\delta}$.
Therefore,  $B({\boldsymbol\delta})$ has rank $\leq n-1$.
But since $\boldsymbol\delta$ defines a birational map, not all $(n-1)$-minors vanish modulo $\det(B^t)$.
Thus,  $B({\boldsymbol\delta})$ has rank $n-1$.
For even more reason, the rank of the Jacobian dual matrix of
${\bf D}$ is $\geq n-1$ (hence $=n-1$, its maximal possible value).
Using again the criterion of \cite{AHA} we derive that ${\bf D}$ defines a Cremona map.

\smallskip

Now, we prove the additional statement of this item.
Let $s$ denote the minimal number of generators of the Rees ideal of ${\bf D}$
of bidegree $(1,*)$, with $*$ representing any value $\geq 1$.
Then  the full Jacobian dual matrix of  ${\bf D}$ is an $s\times n$ matrix over $k[\YY]$ of rank $n-1$
which, as we have just shown, contains the $n\times n$ submatrix $B({\boldsymbol\delta})$.
By \cite[Theorem 2.18, Supplement]{AHA} we know that the inverse map to the Cremona map defined by  ${\bf D}$
takes as its coordinate functions the $(n-1)$-minors of any $(n-1)\times n$ submatrix of rank $n-1$ of
the Jacobian dual matrix of ${\bf D}$, further divided by their $\gcd$. Since $B({\boldsymbol\delta})$
has rank $n-1$, one can take, say, the
submatrix of $B({\boldsymbol\delta})$ formed with the first $n-1$
rows of $B({\boldsymbol\delta})$. Write  $\partial_i(\ZZ)$ for the $(n-1)$th minor
omitting the $i$th column. Then we get $\partial_i(\ZZ) = B_{ni}({\boldsymbol\delta}) =
B_{ni}^t({\boldsymbol\delta}) = X_nd_i(\ZZ)$, where $d_i(\ZZ)$ as before denotes the corresponding source inversion factor
of the birational map defined by $\boldsymbol\delta$.
It follows that $(d_1(\ZZ):\cdots :d_n(\ZZ))$ defines the inverse map to $\mathfrak{D}$.

(iii) By Proposition~\ref{m_greater_than_n_is _birational}, one has $(I^{n-1}, D_1,\ldots,D_n)\subset I^{(n-1)}$.
On the other hand, by Proposition~\ref{local_duality},  $I^{(n-1)}/I^{n-1}$ is minimally generated by $n$
elements of degree $n(n-1)-1$.
To conclude that the  residues of $D_1,\ldots,D_n$ on $I^{(n-1)}/I^{n-1}$
form a set of minimal generators of the latter it suffices to show that they are $k$-linearly independent.
By part (i) they are even $k$-algebraically independent.

(iv) By (i) and (ii), $\{d_1=d_1(\ZZ),\ldots,d_n=d_n(\ZZ)\}$ generate an ideal of codimension $2$
defining the inverse map to  $\mathfrak{D}$.
Write
\begin{equation*}
h_i=h_i(Z_1,\ldots,Z_n):=Z_id_i\;(=B^{t}_{ii}(\boldsymbol\delta)), i=1,\ldots, n
\end{equation*}
Evaluate $h_i$ on $X_i{\bf D}=(X_iD_1,\ldots, X_iD_n)$ (i.e., through $Z_j\mapsto X_iD_j$):
\begin{eqnarray*}
h_i(X_iD_1,\ldots,X_iD_n)&=&X_iD_id_i(X_iD_1,\ldots,X_iD_n)\\
&=&X_i^{n(n-1)}D_id_i(D_1,\ldots,D_n)\\
&=&X_i^{n(n-1)+1}D_i G
\end{eqnarray*}
where $G:=X_i^{-1}\,d_i(D_1,\ldots,D_n)$ is the source inversion factor of the Cremona map defined by $\bf D$.

On the other hand, one has
\begin{eqnarray*}
h_i(X_iD_1,\ldots,X_iD_n)&=&B^{t}_{i,i}(\delta_1(X_i{\bf D}),\ldots,\delta_n(X_i{\bf D}))\\
&=&B^t_{i,i}(E_{i}({\bf\Delta})\Delta_1,\ldots,E_i({\bf\Delta})\Delta_{n+1})\\
&=&E_i({\bf \Delta})^{n-1} B^t_{i,i}(\Delta_1,\ldots,\Delta_n)\\
&=&E_i({\bf \Delta})^{n-1}X_iD_i,
\end{eqnarray*}
where $E_i,\, i=1,\ldots,n$ are a complete set of target inversion factors of the birational map defined by $\boldsymbol\delta$,
as in (\ref{igualdadeestrutural}).
This implies the relation
\begin{equation}\label{shows_degs}
X_i^{n(n-1)}G=E_i({\bf\Delta})^{n-1}.
\end{equation}
Extracting $(n-1)$th roots yields
\begin{equation}\label{fatorinversao}
X_i^{n}G^{1/n-1}=E_i({\bf\Delta})
\end{equation}
Since $E_i$ has degree $n(n-1)-1$ then $(X_1^n,\ldots,X_n^n)G^{1/n-1}\subset I^{n(n-1)-1}$, from which follows
$E:=G^{1/n-1}\in I^{(n(n-1)-1)}.$

\smallskip

The supplementary statement follows from Proposition~\ref{jac_vs_factor} by admitting the irreducibility of $\det(\Theta({\bf D}))$.
This is because as it divides a power of $E$
then it will divide $E$ itself, and since ${\rm deg}(\det(\Theta({\bf D})))={\rm deg}(E)$, they coincide up to a nonzero scalar.
(To hypothetically argue for the irreducibility of $\det(\Theta({\bf D}))$, note that each $D_i$ is an inversion
factor of a Cremona map whose defining coordinates ${\boldsymbol\Delta}$ are sufficiently general forms; for such
a reason one can expect that it too be a sufficiently general polynomial (e.g., because in characteristic zero
it corresponds to a ``general contracted divisor''.) But then also its
partial derivatives are sufficiently general forms, hence $\det(\Theta({\bf D}))$ is an irreducible polynomial,
since having a proper factor is a closed condition on the coefficients and these are products and sums out
of a set of general coefficients.)

(v) We first check that (\ref{resolucaodosDs}) is indeed a complex.
For this, using that $\{D_1,\ldots,D_n\}$ is a complete set of inversion factors of the birational map
defined by $\bf \Delta$, the cofactor matrix of $\Psi$ is
\begin{equation}\label{cofatores}
{\rm adj}(\Psi)=\left(\begin{array}{cccccccccccc}
X_1D_1&X_1D_2&\ldots&X_1D_n\\
X_2D_1&X_2D_2&\ldots&X_2D_n\\
\vdots&\vdots&\ldots&\vdots\\
X_nD_1&X_nD_2&\ldots&X_nD_n
\end{array}\right)
\end{equation}
Since $\Psi$ has rank $n-1$, the cofactor equation gives
\begin{equation}\label{idcofatores1}
{\rm adj}(\Psi)\cdot \Psi=\bf{0}
\end{equation}
and
\begin{equation}\label{idcofatores2}
\Psi\cdot{\rm adj}(\Psi)={\bf 0}
\end{equation}
From (\ref{cofatores}), (\ref{idcofatores1}) implies that $\Psi$ is a matrix of syzygies of $\bf D$, while (\ref{idcofatores2})
gives that $\XX^t$ is a second syzygy thereof.
This shows that one has indeed a complex.
To finish we check the required Fitting codimension by the Buchsbaum--Eisenbud acyclicity criterion.
The verification at the tail of the complex is immediate, while at the middle the codimension of
$$I_{n-1}(\Psi)=I_1({\rm adj}(\Psi))=(\XX)(D_1,\ldots, D_n)$$
is $2$ because (i) showed that the ideal $(D_1,\ldots, D_n)$ has codimension $2$.
\qed

\begin{Remark}\rm Assertion (i) in the last theorem  depends once more on the general linear assumption; thus,
in Example~\ref{symbolic_cyclic} the polarization complex is not exact and, in fact, $\{D_1,D_2,D_3\}$ admit
a proper common factor.
\end{Remark}

\subsubsection{The structure of the symbolic algebra}

Here is the degree numerology so far:

\medskip

$\bullet$ $\deg(d_i)=\deg (D_i)=n(n-1)-1$, for $i=1,\ldots,n$

$\bullet$ $\deg(G)=(n(n-1)-1)n(n-1)-n(n-1)=(n-1)n(n(n-1)-2)$ -- from (\ref{shows_degs}).

$\bullet$ $\deg(E)=\deg(G)/(n-1)=n(n(n-1)-2)$.

\smallskip

%\newpage

Further consideration is given in the following strategic lemma:

\begin{Lemma}\label{gens_of_pi}
Let $\mathcal{L}$ denote an $(n+1)\times n$
general linear matrix over $R=k[X_1,\ldots, X_n]$, with $ n\geq 3$.
Set $I:=I_{n-1}(\mathcal{L})\subset R$
and let $\mathcal{R}^{(I)}$  denote its symbolic Rees algebra.
Let $D_1,\ldots, D_n\in I^{(n-1)}$ and
$E\in I^{(n(n-1)-1)}$ be as above.
Let $\XX=\{X_1,\ldots, X_n\}, \YY=\{Y_1,\ldots, Y_{n+1}\}, \ZZ=\{Z_1,\ldots, Z_n\}, W$ denote  mutually independents
sets of indeterminates.
Consider the surjective homomorphism of $R$-algebras:
$$\pi: k[\XX,\YY,\ZZ,W]\surjects R[It,D_1t^{n-1},\ldots, D_nt^{n-1},Et^{n(n-1)-1}]$$
such that  $X_i\mapsto X_i,$ $Y_j\mapsto \Delta_jt,$ $Z_r\mapsto D_rt^{n-1}$ and $W\mapsto Et^{n(n-1)-1}.$
Then $\ker(\pi)$ contains the following polynomials:
\begin{itemize}
\item[{\rm (1)}] The entries of $\,\XX\cdot B^t$ {\rm (}$n$ such polynomials{\rm )}
\item[{\rm (2)}] The entries of $\,\ZZ\cdot B$ {\rm (}$n$ such polynomials{\rm )}
\item[{\rm (3)}]  The entries of  $\,\XX^t\cdot\ZZ-{\rm adj}(B)$ {\rm (}$n^2$ such polynomials{\rm )}
\item[{\rm (4)}]  The polynomials of the shape $\{X_1W^{n-1}-d_1(\ZZ),\ldots, X_nW^{n-1}-d_n(\ZZ)\}$, where
$d_1,\ldots, d_n$ are forms defining the inverse of $D_1, \ldots, D_n$ {\rm (}$n$ such polynomials{\rm )}
\item[{\rm (5)}]  The polynomials of the shape $\{Y_1W-\delta_1(\ZZ),\ldots,Y_{n+1}W-\delta_{n+1}(\ZZ)\}$, coming from {\rm (\ref{eqdasimb})}
below  {\rm (}$n+1$ such polynomials{\rm )}.
\end{itemize}
\end{Lemma}
%\vspace{-10pt}
\demo
The first four blocks were discussed before, namely:

\noindent (1) These are equations defining the Rees algebra $R[It]$ of $I$ on the polynomial ring $k[\XX,\YY]$.
Since $R[It]$ is a subalgebra of $R[It,D_1t^{n-1},\ldots, D_nt^{n-1},Et^{n(n-1)-1}]$,
then the equations obviously vanish under $\pi$.

\noindent (2)  Note that the matrix $B$ evaluated by $Y_j\mapsto \Delta_jt$
is a syzygy matrix of $\{D_1,\ldots,D_n\}$ by Theorem~\ref{The D's} (iv).
Since $Z_i$ maps to $D_i t^{n-1}$ the vanishing of $I_1(\ZZ\cdot B)$ is clear as well by the same token.

\noindent (3) One argues as in the previous item based on the  proof of  Theorem~\ref{The D's} (iv).

\noindent (4) These equations under $\pi$ just express the fact that $G=E^{n-1}$ is inversion factor of the
Cremona map defined by $\{D_1,\ldots,D_n\}$.

\noindent (5) To discuss these equations, recall the relation obtained in
(\ref{igualdadeestrutural}):
 \begin{eqnarray}
\delta_{j}(B^t_{n1}({\bf \Delta}),\ldots, B^t_{nn}({\bf\Delta}))&=&E_{n}({\bf \Delta})\Delta_{j}
\end{eqnarray}
On the other hand, we have

\begin{eqnarray*}
\delta_j(B^t_{n1}({\bf \Delta}),\ldots, B^t_{nn}({\bf\Delta}))&=&\delta_j(B_{n1}({\bf \Delta}),\ldots, B_{nn}({\bf\Delta}))\\
&=&\delta_j(X_nD_1,\ldots,X_nD_n) \\
&=& X_n^n \delta_j(D_1,\ldots,D_n).
\end{eqnarray*}
Therefore,
\begin{equation*}
E_{n}({\bf \Delta})\Delta_{j}=X_n^n\delta_j({\bf D})
\end{equation*}
Collecting the two resulting expressions yields

\begin{equation*}
X_n^n\Delta_jE=X_n^n\delta_j({\bf D})
\end{equation*}
and hence

\begin{equation}\label{eqdasimb}
\Delta_jE=\delta_j({\bf D})
\end{equation}
as was to be shown.
\qed

\medskip

We note that the intended generator of symbolic order $n(n-1)-1$ is $E$ and not its $(n-1)$th power $G$; this
raises a suspicion as to whether the polynomials of type (4) above are minimal generators of $\ker(\pi)$.
And indeed, we have the following tightening result:

\begin{Proposition}\label{lowering_degs}
Keeping the notation of the previous lemma, in the generation of the ideal $\ker(\pi)$ one may replace
the $n$ equations of the form $X_iW^{n-1}-d_i(\ZZ)$ by $n$ equations of the form $X_iW- Q_i(\YY,\ZZ)$,
where $Q_i(\YY,\ZZ)$ is a polynomial in $k[\YY,\ZZ]$ of the shape
\begin{equation}
Q_i(\YY,\ZZ)=\sum_{
\begin{array}{c}
\{j_1,\ldots, j_{n-2}\}\subset \{1,\ldots,n+1\}\\
t_1+\ldots+t_{n-2}=n-2
\end{array}
}
Y_{j_1}^{t_1}Y_{j_2}^{t_2}\cdots Y_{j_{n-2}}^{t_{n-2}} P_{t_1,\ldots,t_{n-2}}(\ZZ).
\end{equation}
In particular, $(\XX)E\subset I^{n-2}\,(I^{(n-1)})^{n-1}$.
\end{Proposition}
\demo
Let as above $\delta_1=\delta_1(\ZZ),\ldots, \delta_{n+1}=\delta_{n+1}(\ZZ)$ denote the $n$-minors of the matrix $\mathcal{L}'$
and let $\pi$ as be as given.

We claim that for any collection of non-negative integers $ t_1,\ldots,t_s$, with $s\leq n+1$, and for every subset
$\{j_1,\ldots, j_s\}\subset \{1,\ldots, n+1\}$, the polynomials
$$
Y_{j_1}^{t_1}\cdots Y_{j_s}^{t_s} W^{t_1+\ldots+t_s}-\delta_{j_1}(\ZZ)^{t_1}\cdots\delta_{j_s}(\ZZ)^{t_s}\in \ker(\pi)
$$
belong to the ideal generated by the polynomials from block (5) in the statement of the previous lemma.

We proceed by induction on $s$.

The result is clear for $s=1$ because $Y_jW-\delta_j\in \ker(\pi)$ by the previous lemma and is a factor of $Y_j^tW^t-{\delta_j}^t$, for any $t$.

Thus, assume that $s>1$ and that, without loss of generality, $t_1\neq 0$ (the result is trivially satisfied if all $t$'s are null).
Write
{\small
\begin{eqnarray*}
\bigl(Y_{j_1}^{t_1}W^{t_1}-\delta_{j_1}^{t_1}\bigr)&&\kern-20pt Y_{j_2}^{t_2}\cdots Y_{j_{s}}^{t_{s}}W^{t_2+\ldots+t_{s}}=
Y_{j_1}^{t_1}Y_{j_2}^{t_2}\cdots Y_{j_{s}}^{t_{s}} W^{t_1+\ldots+t_s}-
\delta_{j_1}^{t_1}Y_{j_2}^{t_2}\cdots Y_{j_{s}}^{t_{s}}W^{t_2+\ldots+t_{s}}\\
&=&Y_{j_1}^{t_1}Y_{j_2}^{t_2}\cdots Y_{j_{s}}^{t_{s}} W^{t_1+\ldots+t_s}
-\delta_{j_1}^{t_1}\cdots\delta_{j_s}^{t_s}+\delta_{j_1}^{t_1}\cdots\delta_{j_s}^{t_s}
-\delta_{j_1}^{t_1}Y_{j_2}^{t_2}\cdots Y_{j_{s}}^{t_{s}}W^{t_2+\ldots+t_{s}}\\
&=&\bigl(Y_{j_1}^{t_1}Y_{j_2}^{t_2}\cdots Y_{j_{s}}^{t_{s}} W^{t_1+\ldots+t_s}-\delta_{j_1}^{t_1}\cdots\delta_{j_s}^{t_s}\bigr)
-\delta_{j_1}^{t_1}\bigl(Y_{j_2}^{t_2}\cdots Y_{j_{s}}^{t_{s}}W^{t_2+\ldots+t_{s}}-\delta_{j_2}^{t_2}\cdots\delta_{j_s}^{t_s}\bigr)
\end{eqnarray*}
}
Applying the inductive hypothesis on the two ends of this strand of inequalities shows that the polynomial
$$Y_{j_1}^{t_1}Y_{j_2}^{t_2}\cdots Y_{j_{s}}^{t_{s}} W^{t_1+\ldots+t_s}-\delta_{j_1}(\ZZ)^{t_1}\cdots\delta_{j_s}(\ZZ)^{t_s}$$
also belongs to $\ker(\pi)$.
In particular, taking $s=n-2$ and $t_1,\ldots, t_{n-2}$ any partition of $n-2$, the polynomial
\begin{equation}
Y_{j_1}^{t_1}Y_{j_2}^{t_2}\cdots Y_{j_{n-2}}^{t_{n-2}} W^{n-2}-\delta_{j_1}(\ZZ)^{t_1}\cdots\delta_{j_{n-2}}(\ZZ)^{t_{n-2}}
\end{equation}
 belongs to $\ker(\pi)$.

On the other hand, as seen in Theorem~\ref{The D's} (ii), the coordinate forms $\{d_1=d_1(\ZZ),\ldots, d_n=d_n(\ZZ)\}$ defining the inverse
of the Cremona map defined by $\{D_1,\ldots,D_n\}$ also constitute a complete set of source inversion factors
of the birational map defined by the $n$-minors $\delta_1,\ldots,\delta_n$ of the general linear matrix $\mathcal{L}'$.
Therefore, Proposition~\ref{m_greater_than_n_is _birational} gives
\begin{equation*}
(d_1,\ldots, d_n)\subset (\delta_1,\ldots,\delta_{n+1})^{(n-1)}
\end{equation*}
and, for even more reason
 \begin{equation}\label{cont}
 (d_1,\ldots, d_n)\subset (\delta_1,\ldots,\delta_{n+1})^{(n-2)}= (\delta_1,\ldots,\delta_{n+1})^{n-2}.
 \end{equation}

Fixing $i\in \{1,\ldots,n\}$ we can write
 \begin{equation}\label{Q-polynomial}
d_i(\ZZ)=\sum_{t_1+\ldots+t_{n-2}=n-2}P_{t_1,\ldots,t_{n-2}}(\ZZ)\delta_{j_1}(\ZZ)^{t_1}\cdots\delta_{j_{n-2}}(\ZZ)^{t_{n-2}}.
\end{equation}
Thus, one gets that the polynomial
{\scriptsize
\begin{eqnarray*}
X_iW^{n-1}-d_i(\ZZ)&-&\sum_{
t_1+\ldots+t_{n-2}=n-2
}
P_{t_1,\ldots,t_{n-2}}(\ZZ)\,\Bigl(Y_{j_1}^{t_1}Y_{j_2}^{t_2}\cdots Y_{j_{n-2}}^{t_{n-2}} W^{n-2}-\delta_{j_1}(\ZZ)^{t_1}\cdots\delta_{j_{n-2}}(\ZZ)^{t_{n-2}}\,\Bigr)\\[3pt]
&=& W^{n-2}\,\Bigl( X_iW-\sum_{t_1+\ldots+t_{n-2}=n-2}P_{t_1,\ldots,t_{n-2}}(\ZZ)Y_{j_1}^{t_1}Y_{j_2}^{t_2}\cdots Y_{j_{n-2}}^{t_{n-2}}\,
\Bigr),
\end{eqnarray*}
}
for arbitrary subsets $\{j_1,\ldots, j_{n-2}\}\subset \{1,\ldots,n+1\}$,
belongs to $\ker(\pi)$.
Since $\ker(\pi)$ is a prime ideal and $W\not\in \ker(\pi)$, we conclude that
$$X_iW-\sum_{t_1+\ldots+t_{n-2}=n-2}P_{t_1,\ldots,t_{n-2}}(\ZZ)Y_{j_1}^{t_1}Y_{j_2}^{t_2}\cdots Y_{j_{n-2}}^{t_{n-2}}\in \ker(\pi)$$
as was to be shown.

The second statement is clear.
\qed

\medskip

We now come to the main result of this part.

\begin{Theorem}\label{symbolic_m=n+1} {\rm (char}$(k)=0${\rm )}
Let $\mathcal{L}$ denote an $(n+1)\times n$
general linear matrix over $R=k[X_1,\ldots, X_n]$, with $ n\geq 3$.
Set $I:=I_{n-1}(\mathcal{L})\subset R$
and let $\mathcal{R}^{(I)}$  denote its symbolic Rees algebra.
 Let $\pi: R[\YY,\ZZ,W]\surjects R[It,D_1t^{n-1},\ldots, D_nt^{n-1},Et^{n(n-1)-1}]$ stand for the
$R$-algebra homomorphism as defined in {\rm Lemma~\ref{gens_of_pi}}. Then

\begin{enumerate}
\item[{\rm (a)}] The kernel of $\pi$ is the ideal $\mathcal{P}$ generated by the polynomials
{\footnotesize
$$ I_1(\XX\cdot B^t),\; I_1(\ZZ\cdot B),\; I_1(\XX^t\cdot \ZZ- {\rm adj}(B)),\;
 Y_jW-\delta_j(\ZZ)\; (1\leq j\leq n+1),\;  X_iW- Q_i(\YY,\ZZ)\; (1\leq i\leq n),$$
 }
 where $Q_i(\YY,\ZZ)$ is described in {\rm Proposition~\ref{lowering_degs}}.
\item[{\rm (b)}]  $\mathcal{R}^{(I)}=R[It,D_1t^{n-1},\ldots,D_nt^{n-1},Et^{n(n-1)-1}]$
\end{enumerate}
\end{Theorem}
\demo  We first claim that $W$ is a nonzerodivisor on $R[\YY,\ZZ,W]/\mathcal{P}$.
For this, we will use Gr\"obner basis theory.
Namely, consider the degrevlex order with $\ZZ > \YY > \XX > W$.
As is well-known, it suffices to show that $W$ is not a factor of a minimal generator of in$(\mathcal{P})$.
Now, none of the monomials $Y_jW \,(1\leq j\leq n+1),\, X_iW\, (1\leq i\leq n)$ is a minimal generator
of  in$(\mathcal{P})$ since the order first breaks a tie by the degree, while both $\delta_j(\ZZ),\, Q_i(\YY,\ZZ)$
have degree at least $n\geq 3$.
However, a multiple thereof could be a fresh generator of in$(\mathcal{P})$.
We must exclude this possibility.

For subsequent frequent use, we single out the following fact: any $\delta_j=\delta_j(\ZZ)$ is an irreducible
polynomial. For $n\geq 4$ this follows from the fact that $\delta_j$ is a minimal generator of the prime ideal
$I_n(\mathcal{L}')$ (Proposition~\ref{general_cod2} and the general linear nature of $\mathcal{L}'$).
For $n= 3$, the ideal $I_n(\mathcal{L}')$ is only radical, hence one needs a more direct approach.
We may assume that $\mathcal{L}'$ -- just as $\mathcal{L}$ -- being a general linear matrix, up to
sufficiently general elementary row operations, has the form
$$
\left(
\begin{array}{ccc}
0  &   \ell_1     &  \ell_2 \\
\ell_4  &  0    &   \ell_3  \\
\ell_5  &  \ell_6    &  0 \\
X_1   &  X_2    &  X_3
\end{array}
\right),
$$
where the six $\ell_i$'s constitute mutually general $1$-forms.
Since the $\ell_i$'s are general $1$-forms and the minors involving the last row have a similar shape,
it suffices to consider the minors of the first $3$ rows and the one of the last $3$ rows.
These are, respectively:
\begin{equation}\label{3-minors}
\delta_1=X_1\ell_3\ell_6-X_2\ell_3\ell_5+X_3\ell_4\ell_6,\; \delta_4=\ell_1\ell_3\ell_5+\ell_2\ell_4\ell_6.
\end{equation}
%Write $\ell_i=\alpha_{i1}X_1+\alpha_{i2}X_2+\alpha_{i3}X_3$ for $1\leq i\leq 6$, where the eighteen coefficients $\alpha_{ij}$'s are general
%elements of $k$.
%As for $\delta_4$ it is a complete cubic form whose coefficients are
%$$\bigl\{\alpha_{1i}\alpha_{3j}\alpha_{5l}+\alpha_{2i}\alpha_{4j}\alpha_{6l}\,|\, i,j,l\in \{1,2,3\}\,\bigr\}.$$
%Thus, it suffices to argue that these ten coefficients are general as well.
Now replacing every $\ell_i$ by a new variable $Y_i$, the corresponding minors become
$$(X_1Y_6-X_2Y_5)Y_3+X_3Y_4Y_6,\; Y_1Y_3Y_5+Y_2Y_4Y_6,$$
respectively.
The first polynomial is irreducible since it is a primitive polynomial with respect to the variable $Y_3$.
The second is irreducible since it is a binomial whose terms have $\gcd=1$.
Since the $\ell_i$'s are general then mapping $Y_i\mapsto \ell_i$ shows that the polynomials (\ref{3-minors}) are
irreducible as well.

We note {\em en passant} that if the $\ell_i$'s are not general, some minor may have proper factors -- see, e.g,
the matrix (\ref{tcher2}).

\medskip

We now proceed to the Gr\"obner base argument.
Any fresh initial generator is found by an iteration of the so-called $S$-polynomials
(\cite[Section 1.2]{WolmBook2}) associated to pairs of elements of $\mathcal{P}$ starting out with
pairs of the given set of generators thereof.
Since any generator  coming from the part $ I_1(\XX\cdot B^t),\; I_1(\ZZ\cdot B),\; I_1(\XX^t\cdot \ZZ- {\rm adj}(B))$
 does not involve $W$, we must use at least one among the equations $Y_jW \,(1\leq j\leq n+1),\, X_iW\, (1\leq i\leq n)$.
We now analyze the nature of such $S$-polynomials and their iterations stemming from the given starting pair
of generators of $\mathcal{P}$.

\smallskip

%\begin{enumerate}
%\item[{\rm (1)}]

{\bf (1)} Starting pair $\{Y_jW-\delta_j(\ZZ),\, Y_kW-\delta_k(\ZZ)\}$ $(j\neq k)$

\smallskip

Consider the respective initial terms, which are pure monomials in $\ZZ$
 -- this is because, as already remarked, $\deg(\delta_j(\ZZ))= n>2$.
 Say, $M_j=M_j(\ZZ), M_k=M_k(\ZZ)$ are the respective initial terms and set $H:=\gcd(M_j,M_k)$,
 so $M_j=N_jH$, $M_k=N_kH$, with $\gcd(N_j,N_k)=1$.
Then the associated $S$-polynomial has the shape
\begin{equation}\label{Spoly}
S:=(N_kY_j-N_jY_k)W - (N_k\delta_j'-N_j\delta_k'),
\end{equation}
where $\delta_j=M_j+\delta_j', \delta_k=M_k+\delta_k'$.
By a similar token, $\deg(N_k\delta_j')=\deg(N_j\delta_k')>\deg(N_kY_jW)=\deg(N_jY_kW)$, and hence
the initial term of (\ref{Spoly}) has to come from either $N_k\delta_j'$ or $N_j\delta_k'$, provided
we make sure that the  $N_k\delta_j'-N_j\delta_k'$ does not vanish.
But since $N_jM_k-N_kM_j=0$  by construction, this vanishing would imply the relation $N_k\delta_j-N_j\delta_k=0$.
 However,  $\delta_j, \delta_k$ are non-associate irreducible polynomials,
hence are relatively prime. This would force a trivial relation, hence $N_j$ would be
multiples of $\delta_j$ -- this is absurd since $N_j$ is a monomial.

Repeat the $S$-polynomial procedure using (\ref{Spoly}) and any other equation of type $Y_pW-\delta_p$
obtaining a new $S$-polynomial.
To make it explicit, say, the initial term of old $S$ comes from $N_k\delta_j'$; then write $\delta_j'=M_j'+\delta_j''$,
where $M_j'$ is the initial term.
Also write, as above, $\delta_p=M_p+\delta_p'$, with $M_p$ its initial term.
Finally, set $M_j'=N_j'H,\, M_p=N_p'H$, where $\gcd(N_j',N'_p)=1$.
Then the updated $S$-polynomial is
\begin{eqnarray*}
S':= \Bigl(N_p'(N_k Y_j-N_j Y_k)- N_kN_j'Y_p\Bigr)W\\
-\Bigl (N_p'(N_k\delta_j'')-N_p'(N_j\delta_k')-N_k(N_j'\delta_p')\Bigr).
\end{eqnarray*}
Counting degrees as before, we that the degree of the top part is lower than that of the bottom part.
Therefore, the initial term of $S'$ will come off the bottom part provided we show it does not vanish.
Supposing this were the case, using the basic  $S$-pair relation $N_p'(N_kM_j')-(N_kN_j')M_p$, we get
the relation $N_p'(N_k\delta_j'-N_j\delta_k')-N_kN_j'\delta_p=0$.
Since $\gcd(N_j',N'_p)=1$, we see that $N'p$ divides $N_k\delta_p$. But $\delta_p$ is an irreducible, hence $N'_p$ divides $N_k$.
Substituting in the previous relation and simplifying yields $N_k\delta_j'-N_j\delta_k'=N\delta_p$,
for some monomial $N\in k[\ZZ]$.
But this implies that our initial $S$-polynomial in (\ref{Spoly}) has the form $(N_kY_j-N_jY_k)W - N\delta_p$.
Using $Y_p-\delta_p$, one gets $(N_kY_j-N_jY_k-NY_p)W=0$ and hence $N_kY_j-N_jY_k-NY_p=0$.
This is nonsense since $N_k,N_j,N$ are polynomials in $k[\ZZ]$.

\smallskip

Now the general iterated step is clear, therefore the iteration of $S$-polynomials using only this packet
of equations gives fresh initial generators which are
monomials in $\ZZ$ exclusively.

\medskip

{\bf (2)} Starting pair $\{X_iW- Q_i(\YY,\ZZ),\, X_lW- Q_l(\YY,\ZZ)\}$ $(i\neq l)$

\smallskip

The argumentative strategy is analogous to the one in the previous case: write
$$\left\{
\begin{array}{cc}
Q_i=M_i+Q_i', & Q_l=M_l+Q_l'\\
M_i=M_i'H, & M_l=M_l'H,
\end{array}
\right.
$$
where $M_i={\rm in}(Q_i),\, M_l={\rm in}(Q_l)$ and $\gcd(M_i',M_l')=1$.
Note that, from (\ref{Q-polynomial}) and since $\deg(d_i(\ZZ))=n(n-1)-1$ and $\deg(\delta(\ZZ))=n(n-2)$,
one has $\deg(Q_i)=2n-3$.
Then the resulting polynomial is the sum of two homogeneous polynomials
$$S:=(M_l'X_i-M_i'X_l)W-(M_l'Q_i'-M_i'Q_l'),$$
where $\deg (M_l'X_iW)=2n-3-h+2=2n-1-h,\, \deg(M_l'Q_i')=2n-3-h+2n-3=4n-6-h$, with $h=\deg(H)$.
Again, we have the strict inequality $4n-6> 2n-1$, for $n\geq 3$.
To show that the initial term of the above polynomial belongs to the rightmost polynomial we need to
know that the latter does not vanish.
But if it did, then we would have the equality $M_l'Q_i=M_i'Q_l$, where $\gcd(M_i',M_l')=1$.
Now, since the multipliers $M_i',M'_l$ are relatively prime then $M_i'$ is a factor of $Q_i$.
By (\ref{Q-polynomial}), evaluating  we would get that $d_i$ has a monomial factor in $k[\ZZ]$; this is ruled out by
fact that $d_i$ is a coordinate function of the inverse map to the Cremona map defined by $\{{D}_1,\ldots {D}_n\}$
which are sufficiently general forms.

\smallskip

We can now iterate as in case (1).
Thus, let
\begin{equation}\label{S_X}
S_X=\Bigl(\sum_i N_i(\YY,\ZZ)X_i\Bigr)W-P(\YY,\ZZ)
\end{equation}
stand for an iterated $S$-polynomial out of the ``$X_iW$'' packet, with $M=M(\YY,\ZZ)$ denoting the corresponding
initial term.
By induction, we have $\deg(N_i)+2 <\deg(P)$.
Write $P=M+P'=:{\rm in}(P)+P'$ and $Q_r=M_r+Q_r':={\rm in}(Q_r)+Q_r'$. Then the new $S$-polynomial has the shape
$$\sum_i M'_rN_iX_iW - M'X_rW - (M_r'P'-M'Q_r'),$$
where $M=M'H,\, M_r=M'_rH$, with $\gcd(M',M'_r)=1$.
We assume that $\deg (H)>0$ as otherwise there is nothing to prove by \cite[Exercise 1.2.2]{WolmBook2}.
Then $\deg(M'_rN_iX_iW)=2n-3+\deg(N_i)+2 -h< 2n-2+\deg(P)-h=\deg(M_r'P')$ and, similarly, $\deg(M'X_rW)=\deg(M)+2-h=
\deg(P)+2-h<\deg(P)+2n-3-h=\deg(M_r'P')$.
Moreover, if $M_r'P'=M'Q_r'$ then $M_r'P=M'Q_r$ as well.
If $\gcd(P,Q_r)=1$ then $M_r'$ must be a multiple of $Q_r$, which is impossible since $\deg(M'_r)< \deg(Q_r)$ by
hypothesis.
Then $P$ and $Q_r$ must have a proper common factor.
Now, since the multipliers $M_r',M'$ are relatively prime then $M_r'$ is a factor of $Q_r$.
Under $Y_j\mapsto \delta_j(\ZZ)$  we would get that $d_r$ has a monomial factor in $k[\ZZ]$; this is again ruled out
as above.

Thus, the initial term of an $S$-polynomial from pairs consisting of any previous $S$-polynomial obtained and any other
equation $X_rW- Q_r(\YY,\ZZ)$ is a monomial in $\YY$ and $\ZZ$ alone.

\medskip

{\bf (3)}  (Mixed starting pair) One of the pairs
$$\{Y_jW-\delta_j(\ZZ),\,S_X\}\quad {\rm or} \quad \{ X_iW- Q_i(\YY,\ZZ),\, S_Y\},$$
for some $1\leq j\leq n+1$ and some $1\leq i\leq n$, where $S_Y$ (respectively, $S_X$) is any $S$-polynomial from the
``$YW$'' packet (respectively, from the ``$XW$'' packet).

Let us deal with these pairs separately.
For the first pair, let $S_X$ have the expression as in (\ref{S_X}).
Then the new $S$-polynomial has the form
$$\sum_i M'_jN_iX_iW - M'Y_jW - (M_j'P'-M'\delta_j'),$$
where $\gcd(M',M'_j)=1$.
Degree counting gives $\deg(M'_jN_iX_iW )=n+\deg(N_i)+2-h<n+\deg(P)-h$.
Moreover, vanishing of the rightmost polynomials would lead to $M_j'P=M'\delta_j$.
As before, we are forced to conclude that $\delta_j$ has a factor which is a monomial.
But this is impossible since $\delta_j$ is irreducible.

\smallskip

For the pair of the second kind, let
$$S_Y=(\sum_j N_jY_j)W-\sum_jN_j\delta_j^{<s_j>}$$
denote an $S$-polynomial as iterated from the ``$Y_jW$'' packet.
Here $\delta_j^{<s_j>}$ denotes a suitable summand of $\delta_j$ and the initial term of $S_Y$.
Form the $S$-polynomial with some $X_iW-Q_i$, $Q_i=Q_i(\YY,\ZZ)$, getting:
$$M_i'\bigl(\sum_j N_jY_j\bigr)W-M'_{j_0}X_iW -\Bigl(M_i'\sum_j N_j\delta_j^{<s'_j>}-M'_{j_0}Q'_i\Bigr),$$
where

$$
\left\{\begin{array}{rl}
M_{j_0}= & {\rm in}(\delta_j^{<s_{j_0}>})\\
 M_i= & {\rm in}(Q_i)
 \end{array}
 \right.
 \quad \mbox{\rm and}\quad
 \left\{\begin{array}{rc}
  N_{j_0}M_{j_0}=&M'_{j_0}H\\
   M_i=&M'_iH
\end{array}
\right.
$$
with
$\gcd(M'_{j_0}, M'_i)=1$ and $\delta_j^{<s'_j>}$ are the updated summands of $\delta_j$.
Once again, an immediate degree count tells us that the initial term of the new $S$-polynomial
comes from the right most difference above, as long as the latter does not vanish.
Supposing it did,  we would as before get the relation $M_i'\sum_j N_j\delta_j^{<s_j>}=M'_{j_0}Q_i$,
with monomial multipliers relatively prime.
This implies that $M'_i$ is a factor of $Q_i$, which leads to a relation $\sum_j N_j\delta_j^{<s_j>}=M'_{j_0}Q''_i$
with, say, $Q_i=M_i Q''_i$.
Substituting back in $S_Y$ gives $S_Y=(\sum_j N_jY_j)W-M'_{j_0}Q''_i$.
Multiplying $S_Y$ by $M_i$ and $X_iW-Q_i$ by $M'_{j_0}$, and subtracting yields $(\sum_j M'_iN_jY_j-M'_{j_0}X_i)W=0$.
Therefore, $\sum_j M'_iN_jY_j=M'_{j_0}X_i$, which implies that $M'_{j_0}$ belong to the ideal generated by a nonempty
subset of the of the $\YY$ variables; this is absurd since $M'_{j_0}\in k[\ZZ]$.

To conclude these cases,  note that iterating these two types of $S$-polynomials, we obtain similarly that any pair $\{S_Y,S_X\}$
yields an $S$-polynomial whose initial term is not divisible by $W$.

\medskip

{\bf (4)}  Starting pair $\{Y_jW-\delta_j(\ZZ),\, q\}$

\smallskip

Here $q$ is a generator out of $ I_1(\XX\cdot B^t),\; I_1(\ZZ\cdot B),\; I_1(\XX^t\cdot \ZZ- {\rm adj}(B))$.

Let $q$ come from $ I_1(\XX\cdot B^t)$. Then its initial term is of the form $\alpha Y_kX_i$.
Since the initial term of $Y_jW-\delta_j(\ZZ)$ is a monomial in $\ZZ$ alone, these two monomials are relatively
prime.
Therefore, the resulting $S$-polynomial reduces to zero relative to the pair $\{Y_jW-\delta_j(\ZZ),\, q\}$
(\cite[Exercise 1.2.2]{WolmBook2}) and hence, produces no fresh initial generator.

Assume now that $q$ comes from the packet $I_1(\ZZ\cdot B)$. By a similar token, the initial term of $q$ has the
form $\beta Z_iY_k$.
Since the initial term of $Y_jW-\delta_j(\ZZ)$ is a monomial in $\ZZ$ alone, the only way to get away from reducing
to zero as before is that $\beta Z_i$ divide this $\ZZ$-monomial.
Thus, let $ M(\ZZ)Z_i$ denote the initial term of $\delta_j(\ZZ)$ and write $\delta_j=\beta M(\ZZ)Z_i+P(\ZZ)$.
Then the resulting polynomial is
$$S:= Y_kY_jW-Y_kP(\ZZ)+M(\ZZ) q(\ZZ,\YY),$$
where $q(\ZZ,\YY)$ is a $2$-form of bidegree $(1,1)$ in $\ZZ,\YY$.
Clearly, $3=\deg(Y_kY_jW)<1+n=\deg(Y_kP(\ZZ))=\deg(M(\ZZ) q(\ZZ,\YY))$, hence the initial term of $S$ involves
only $\ZZ$ and $\YY$ variables provided we show that $-Y_kP(\ZZ)+M(\ZZ) q(\ZZ,\YY)$ does not vanish.
Now, a similar reasoning as employed at the end of the argument of (1), shows that this vanishing entails
a  monomial syzygy between $\delta_j(\ZZ)$ and the quadric $q$ with relatively prime multipliers.
This then forces $\delta_j(\ZZ)$ and $q$ to have a common factor. But $q$ is bihomogenous, so a common factor
would have to be a variable $Z_l$. On the other hand, $\delta_j(\ZZ)$ is irreducible, so
cannot admit such a factor.

Keeping the essential shape of the $S$-polynomial obtained, namely, $S=Y_kY_j-P(\YY,\ZZ)$, with $P(\YY,\ZZ)$ homogeneous
of degree $n+1$ involving effectively both $\YY$ and $\ZZ$ variables, let us iterate with the pair $\{S, Y_rW-\delta_r(\ZZ)\}$,
for given $1\leq r\leq n+1$.
Write $P(\YY,\ZZ)=M(\YY,\ZZ)+P'(\YY,\ZZ),\, \delta_r(\ZZ)=N(\ZZ)+\delta_r'(\ZZ)$, where $M=M(\YY,\ZZ), N=N(\ZZ)$ are the respective
initial terms, and $M=M'H, N=N'H$, with $\gcd(M',N')=1$.
Then the new $S$-polynomial is
$$ \bigl(N'Y_kY_j-M'Y_r)W - (N'P'-M'\delta_r'),
$$
where the leftmost polynomial is homogeneous of degree $n+3-h$, with $h=\deg(H)$, while the rightmost polynomial has
degree $2n+1-h> n+3-h$, for $n\geq 3$.
On the other hand, the rightmost polynomial is nonzero because, otherwise, it would imply that $N'P=M'\delta_r$.
Since $N',M'$ are relatively prime, $\delta_r$ would be a multiple of $P=P(\YY,\ZZ)$.
But this is absurd since $\delta_r\in k[\ZZ]$ while $P\notin k[\ZZ]$.
Therefore, the initial term of the updated $S$-polynomial comes from the rightmost polynomial and does not involve $W$.
The inductive procedure is now clear: the ``new'' $S$-polynomial is a sum of two polynomials, the first involving $W$ and
degree growing like $(s-1)n+3-t$, for $s\geq 2$ and some $t\geq 0$, the second a nonzero polynomial involving effectively
the variables $\YY, \ZZ$ and with degree
growing like $sn+1-t > (s-1)n+3-t$ (for $n\geq 3$).

\smallskip

Finally, consider the case where $q$ comes from the packet $I_1(\XX^t\cdot \ZZ- {\rm adj}(B))$.
If $n\geq 4$, the initial term is decided by degree and has to come from some cofactor of $B$ -- the latter
having degree $n-1\geq 3>2=\deg(X_iZ_l)$, for any choice of $i,l$.
In this case, once again, the $S$-polynomial reduces to zero.
Finally, let $n=3$. Since we are assuming the revlex order upon monomials of same degree,  the initial term
of $P$ comes from a cofactor of $B$, so we are done again.

\begin {Remark}\label{the_end}\rm
To close this case, we ought to consider the $S$-polynomial from the pair consisting of a polynomial
of the ``$Y_jW$'' packet and some previous $S$-polynomial among one of the three kinds.
But, as we have seen, the only iterated $S$-polynomials that play any role come from the pairs
$$\{Y_jW-\delta_j(\ZZ),\, q\in I_1(\ZZ\cdot B)\}.$$
One can see that this iteration follows a pattern analogous to the first iterate, in which the initial term
lives in $k[\YY,\ZZ]$.
\end{Remark}

\medskip

{\bf (5)}  Starting pair $\{X_iW-Q_i(\YY,\ZZ),\, q\}$

Here $q$ is again  a generator out of $ I_1(\XX\cdot B^t),\; I_1(\ZZ\cdot B),\; I_1(\XX^t\cdot \ZZ- {\rm adj}(B))$.

\smallskip

The initial term of $X_iW-Q_i(\YY,\ZZ)$ involves both $\YY$ and $\ZZ$.
This breaks the symmetry with respect to the discussion in case (4).

Let first $q$ come from  $ I_1(\XX\cdot B^t)$.
Then $q=\alpha Y_jX_l+q'$, where ${\rm in}(q)=\alpha Y_jX_l$.
Note that $q'\neq 0$ -- i.e., $q$ is not a monomial -- since the entries of a column of $B^t$ are partial $\XX$-derivatives of
minors of a general linear matrix.
Write as before $Q_i=M_i+Q'_i$, where ${\rm in}(Q_i)=M_i= N_i\cdot\alpha Y_j$.
Clearly, $X_l$ does not divide $N_i$.
The resulting $S$-polynomial is $X_lX_iW - (X_lQ'_i-N_iq')$.
One has $\deg(X_lQ'_i-N_iq')=1+2n-3=2n-3-1+2=2n-2>3$, for $n\geq 3$.
Moreover, if $X_lQ'_i=N_iq'$ then $X_lQ_i=N_iq$. This forces $X_l$ to be a factor of $q$, which implies
that $q$ is monomial, contradicting its nature as pointed out.

Now assume that $q$ come from  $ I_1(\ZZ\cdot B)$.
Then $q=\beta Z_kY_j+q'$, where ${\rm in}(q)=\beta Z_kY_j$.
The same remarks about the nature of $q$ hold as above.
Keeping the same notation, $Q_i=M_i+Q'_i$, where ${\rm in}(Q_i)=M_i$.
If $Z_kY_j$ divides $M_i$ altogether, then the resulting polynomial is of the form
$X_iW-(Q'_i-P_iq')$, for suitable $P_i\in k[\YY,\ZZ]$ homogeneous of degree $2n-3-2+2=2n-3>2$.
Thus, we may assume that either $Z_k$ divides $M_i$ and $Y_j$ does not divide $M_i$, or vice versa.
Although  the role of $\YY$ and $\ZZ$ are not quite
 symmetric in the data, the pattern is pretty much the same (and much the same as the previous case).
Say, $M_i=N_i\cdot \beta Z_k$, with $Y_j$ not dividing $M_i$.
The resulting $S$-polynomial is $Y_jX_i-(Y_jQ'_i-N_iq')$.
Again the inequality $2n-2>3$ says that the initial term is part of $Y_jQ'_i-N_iq'$.
Moreover, $Y_jQ_i=N_iq$ would imply that $Y_j$ divide $q$, again a contradiction.

Finally, we settle the last case where  $q$ comes from the packet $I_1(\XX^t\cdot \ZZ- {\rm adj}(B))$.
If $n\geq 4$, the initial term is decided by degree and has to come from some cofactor of $B$ -- the latter
having degree $n-1\geq 3>2=\deg(X_iZ_l)$, for any $i,l$.
Say, $q=C(\YY)+q'$, with ${\rm in}(q)=C=C(\YY)$ of degree $n-1$.
As before, $Q_i=M_i+Q'_i$, where ${\rm in}(Q_i)=M_i$.
Set $C=C'H, M_i=N_i H$, $\gcd(C',N_i)=1$.
The resulting $S$-polynomial is $C'X_iW-(C'Q'_i-N_iq')$, where $\deg (C'Q'_i-N_iq')=n-1+2n-3=3n-4>n+1=\deg(C'X_iW)$.
Furthermore, if $C'Q_i=N_iq$ would imply that $C'$ divide $q$; this is absurd since $q$ is of the form $X_iZ_l-p(\YY)$.

At last, let $n=3$. Since we are assuming the revlex order upon monomials of same degree,  the initial term
of $P$ comes from a cofactor of $B$, so we are done again.

\smallskip

To close this item, we refer to Remark~\ref{the_end}, noting that here one has to consider
iterated $S$-polynomials form all three kinds, as none reduces to zero right at the outset.

\medskip

(a) Let $\mathcal{P}\subset R[\YY,\ZZ,W]$ denote the ideal generated by those many equations in the statement.
By Lemma~\ref{gens_of_pi} and Proposition~\ref{lowering_degs}, we have $\mathcal{P}\subset \ker(\pi)$.
The two ideals have same codimension: $2n+1$. Indeed, the algebra $A:=R[It,D_1t^{n-1},\ldots,D_nt^{n-1},Et^{n(n-1)-1}]$
has the same dimension as the Rees algebra $R[It]$, which is $n+1$; this shows that $\ker(\pi)$ has codimension
$2n+1$.
As for  $\mathcal{P}$, we localize at the powers of $W$.
Then $\mathcal{P}$ and $\mathcal{P}[W^{-1}]\subset k[\XX,\YY,\ZZ, W, W^{-1}]$ have the same codimension.
But in the latter the generators
$$\{Y_j-W^{-1}\delta_j(\ZZ),\, X_i-W^{-1}Q_i(\YY,\ZZ)\,|\, 1\leq j\leq n+1, 1\leq i\leq n\}$$
form a regular sequence of length $n+1+n=2n+1$.

Therefore, to show that $\mathcal{P}= \ker(\pi)$ it suffices to prove that $\mathcal{P}$ is a prime ideal.
By localizing at the powers of $W$, one gets an isomorphism of $k$-algebras
\begin{equation}\label{isomorphism}
k[\XX,\YY,\ZZ, W, W^{-1}]/\mathcal{P}[W^{-1}]\simeq k[\ZZ,W, W^{-1}]/\widetilde{\mathcal{P}[W^{-1}]}
\end{equation}
by mapping $X_i\mapsto W^{-1}Q_i(\YY,\ZZ)$ and subsequently $Y_j\mapsto W^{-1}\delta_j(\ZZ)$.
Since $k[\ZZ,W, W^{-1}]$ has dimension $n+1$, we must conclude that $\widetilde{\mathcal{P}[W^{-1}]}=0$.
(As a control of quality one has that, e.g., $I_1(\XX\cdot B^t)$ maps to
$I_1(d_1(\ZZ)\cdots d_n(\ZZ))\cdot B^t(\delta_1(\ZZ),\ldots, \delta_{n+1}(\ZZ))$, which vanishes
as seen in the proof of Theorem~\ref{The D's} (ii).)

Therefore, $\mathcal{P}[W^{-1}]$ is a prime ideal, and hence so is $\mathcal{P}$.

\smallskip

(b)
We apply to the algebra $A$ the criterion of Vasconcelos (\cite[Propositions 7.1.4 and 10.5.1]{Wolmbook})
mentioned in the first section of this paper (Proposition~\ref{idealtransform}).
By Proposition~\ref{alma_mater} (ii)
the required hypothesis is satisfied -- note the need for characteristic zero at this point.
Therefore, it suffices to prove that the grade of the extended ideal $(\XX)A$
is at least $2$.
For this we claim that the grade of $(\XX)A$ is the same as the grade of its extension to the localization
$A_w$ at the powers of the image $w$ of $W$.
To see this it is enough to show that $w$ avoids some associated prime $\wp$ of $A/(\XX)A$
such that ${\rm grade}((\XX)A)={\rm grade}(\wp)$.
We show more, namely:

{\sc Claim:} $w$ is regular on $A/(\XX)A$.

{\sc Proof.} Now, by (a) we know that $\ker(\pi)=\mathcal{P}$.
Thus, we have to show that $W$ is a non-zero-divisor modulo the larger ideal
\begin{equation*}
\bigl(\XX,\mathcal{P}\bigr)=\bigl(\XX,\; I_1(\ZZ\cdot B),\; I_{n-1}(B),\; Q_i(\YY,\ZZ)\; \mbox{\small $(1\leq i\leq n)$},\;
 Y_jW-\delta_j(\ZZ)\; \mbox{\small $(1\leq j\leq n+1)$}\bigr).
 \end{equation*}
We follow the same Gr\"obner basis line of argument as before to show that $W$ does not divide any generator of
${\rm in}(\XX,\mathcal{P})$.
Note that we can overlook the pairs $\{Y_jW-\delta_j(\ZZ),\, q\}$, with $q\in  I_{n-1}(B)$.
Indeed, this is clear since the initial degree of the first of these polynomials comes from
$\delta_j(\ZZ)$ due to its degree being $n\geq 3$.
The only remaining relevant pairs are then $\{Y_jW-\delta_j(\ZZ),\, q\}$, with $q\in I_1(\ZZ\cdot B)$, and
$\{Y_jW-\delta_j(\ZZ),\, Q_i(\YY,\ZZ)\}$.
The first kind, as well as its descendants, have been dealt with in part (4) and Remark~\ref{the_end}.

Consider the pair $\{Y_jW-\delta_j(\ZZ),\, Q_i(\YY,\ZZ)\}$.
 Say, $M_j=M_j(\ZZ), M_i=M_i(\ZZ)$ are the respective initial terms of $\delta_j=\delta_j(\ZZ)$ and $Q_i=Q_i(\YY,\ZZ)$;
 set $H:=\gcd(M_j,M_i)$,
 so $M_j=N_jH$, $M_i=N_iH$, with $\gcd(N_j,N_i)=1$.
Then the associated $S$-polynomial has the shape
$S:=N_iY_jW - (N_k\delta_j'-N_jQ_i'),
$
where $\delta_j=M_j+\delta_j', Q_i=M_i+Q_i'$.
By a similar token, $\deg(N_i\delta_j')=\deg(N_jQ_i')>\deg(N_iY_jW)$, and hence
the initial term of $S$ has to come from either $N_i\delta_j'$ or $N_jQ_i'$, provided
we make sure that the  $N_i\delta_j'-N_jQ_i'$ does not vanish.
But since $N_jM_i-N_iM_j=0$  by construction, this vanishing would imply the relation $N_i\delta_j=N_jQ_i$.
Evaluating $Y_l\mapsto \delta_l$ for all $l$ and using that $\delta_j$ is irreducible implies that
$\delta_j$ divide $d_i$.
But this is impossible because, as was already remarked,  $d_i$ is a general form.

This analysis is repeated with a pair $\{(N_iY_jW - (N_k\delta_j'-N_jQ_i',\, Q_i(\YY,\ZZ)\}$ and so forth,
by an obvious recursion.
So much for the proof that $W$ is a non-zero-divisor on $A/(\XX)A$.

\smallskip

We are now left with computing ${\rm grade}((\XX)A_w)$.
By (\ref{isomorphism}), $A_w\simeq k[\ZZ,W, W^{-1}]$, which is a Cohen--Macaulay graded ring.
Since the image of $(\XX)A_w$ is a graded ideal, its grade coincides with its codimension.
Now, the image of $(\XX)A_w$ by the isomorphism is the ideal
$$ \bigl(Q_1(W^{-1} \delta_j(\ZZ),\ZZ),\, \ldots, \,Q_n(W^{-1} \delta_j(\ZZ),\ZZ)\,\bigr)\subset k[\ZZ,W, W^{-1}].$$
By homogeneity we can pull out $W^{-1}$. Then (\ref{Q-polynomial}) shows that this ideal is generated
by the coordinate functions $\{d_1,\ldots,d_n\}$ of the inverse map of the Cremona map defined by
$D_1,\ldots,D_{n}$.
By construction, these forms have trivial $\gcd$, hence the ideal they generate has indeed codimension at least $2$.
\qed
\begin{Remark}\label{erratic}\rm
(1) Note that in the proof of part (b) of the theorem it would suffice to show that $A$ satisfies the
Serre property $(S_2)$.
As a matter of fact, one wonders if the symbolic algebra is Cohen--Macaulay, in which case it would be
a Gorenstein normal domain by \cite{ST}.
We have verified this in the case $n=3$ by writing an explicit regular sequence of length $n+1=4$.
The terms of the sequence can actually be taken to be linear forms involving only the $X$ and $Z$ variables
and $W$ -- this exploits the fact that in this dimension one can change to a grading where the $X,Y,Z$ part is standard
and the variable $W$ has weight $2$.
In this grading the Hilbert series is $(1+7t+13t^2+7t^3+t^4)/(1-t)^4$.
For $n\geq 4$ some of such facilitating features are not available.
On the other hand, even for $n=3$, the property that $W$ is a non-zero-divisor on $A/(\XX)A$ is really
on the edge as the ideal $(\XX,\YY,\ZZ)A$ is an associated prime ideal of $A/(\XX)A$.

(2)
In the case $m\geq n+2$ it may happen that elements of $I^{(n-1)}$  have standard degree less than $(m-1)(n-1)-1$.
The simplest such situation occurs with $n=3$ and $m=7$, in which case $I^{(2)}$ admits $3$ minimal generators
of degree $10$.
This implies that, in this range, the inclusion $(\XX)I^{(r)}\subset I^{r}$ for every $r\geq 0$ fails.
This is an indication that,  for general values of $m,n$, it may be difficult to guess bounds for the value of the saturation exponent,
so as to have {\rm Proposition~\ref{alma_mater} (ii)} become more precise.

(3)
Computational evidence showed that in the smallest possible numerology ($n=3,m=5$) the behavior of the symbolic powers
is quite erratic: in the range $2\leq r\leq 5$ there are genuine generators in $I^{(r)}$.
The subsequent symbolic powers have an unpredictable behavior with genuine generators creeping up on irregular
intervals; we found new symbolic generators even in $I^{(23)}$.
It seems reasonable to wonder whether for $m>n+1\geq 4$ the symbolic Rees algebra
$\mathcal{R}^{(I)}$ of $I$ is finitely generated.
\end{Remark}

We close with a couple of more general questions.

\begin{Question}\rm
(1)  It would be interesting to describe classes of (characteristic free) perfect, codimension $2$,  homogeneous, prime ideals
$I\subset R=k[X_1,\ldots,X_n]$, generated in fixed degree  such that
$R/I$ is normal and $I$ admits non-ordinary symbolic powers $I^{(m)}$ of order $m\leq n-2$.

(2) Note that in the setup of this paper, the two alternatives in \cite[Proposition 3.5.13]{WolmBook2}
coincide set-theoretically, namely, the radical of the Jacobian ideal is read off the free presentation of the ideal.
This phenomenon played a central role in the preliminaries of this work.
It seems appropriate to ask when this is the case beyond the present assumptions.
\end{Question}

%%%%%%%%%%%%%%%%%%%%%%%%% bibliografia %%%%%%%%%%%%%%%%%%%%%%%%%%%%

\noindent {\bf Authors' addresses:}

\medskip

\noindent {\sc Zaqueu Ramos}, Departamento de Matem\'atica, CCET, Universidade Federal de Sergipe\\
49100-000 S\~ao Cristov\~ao, Sergipe, Brazil\\
{\em e-mail}: zaqueu@gmail.com\\

\noindent {\sc Aron Simis},  Departamento de Matem\'atica, CCEN, Universidade Federal
de Pernambuco\\
 50740-560 Recife, PE, Brazil.\\
{\em e-mail}:  aron@dmat.ufpe.br

\end{document}